\documentclass[aos]{imsart}

\RequirePackage{amsthm,amsmath,amsfonts,amssymb}
\RequirePackage[numbers]{natbib}
\RequirePackage[colorlinks,citecolor=blue,urlcolor=blue]{hyperref}
\RequirePackage{graphicx}
\usepackage{multirow}
\usepackage[ruled,vlined]{algorithm2e}
\usepackage{tikz}

\usepackage{xr}
\makeatletter
\newcommand*{\addFileDependency}[1]{
  \typeout{(#1)}
  \@addtofilelist{#1}
  \IfFileExists{#1}{}{\typeout{No file #1.}}
}
\makeatother

\newcommand*{\myexternaldocument}[1]{%
    \externaldocument{#1}%
    \addFileDependency{#1.tex}%
    \addFileDependency{#1.aux}%
}

\myexternaldocument{Supplementary_v3}

\startlocaldefs

\newtheorem{theorem}{Theorem}[section]

\newtheorem{cor}[theorem]{Corollary}

\theoremstyle{remark}

\newcommand{\E}{\mathrm{E}}

\newcommand{\am}{{\alpha m}}

\newcommand{\var}{\mathrm{Var}}
\newcommand{\OL}{{\scriptscriptstyle OL}}
\newcommand{\NOL}{\scriptscriptstyle NOL}
\newcommand{\cov}{\mathrm{Cov}}
\newtheorem{Lemma}{Lemma}

\endlocaldefs

\begin{document}

\begin{frontmatter}
\title{On optimal  block resampling  for Gaussian-subordinated long-range dependent processes}

\begin{aug}
\author[A]{\fnms{Qihao} \snm{Zhang}\ead[label=e1]{qihaoz@iastate.edu}},
\author[B]{\fnms{Soumendra N.} \snm{ Lahiri}\ead[label=e2]{s.lahiri@wustl.edu}}
\and
\author[A]{\fnms{Daniel J.} \snm{Nordman}\ead[label=e3]{dnordman@iastate.edu}}
\address[A]{Department of Statistics,
Iowa State University,
\printead{e1,e3}}

\address[B]{Department of Mathematics and Statistics,
Washington University in St. Louis,
\printead{e2}}
\end{aug}

\begin{abstract}
 Block-based resampling estimators have been intensively investigated for weakly dependent
time processes,
which  has helped to inform   implementation (e.g., best block sizes). However,  little is known
about resampling performance and block sizes under strong  or long-range   dependence. To establish guideposts in block selection, we consider   a broad class of strongly dependent
 time processes, formed by a transformation
of a stationary long-memory Gaussian series, and examine block-based resampling estimators for the variance
 of the prototypical sample mean; extensions to general statistical functionals are also considered.   Unlike weak dependence,  the properties of  resampling estimators under strong dependence  are shown to depend intricately on  the nature of non-linearity in the time series (beyond Hermite ranks) in addition the long-memory coefficient and    block size.   Additionally,  the intuition  has often been that optimal block sizes  should be larger under strong dependence (say $O(n^{1/2})$ for a sample size $n$)
 than the optimal order   $O(n^{1/3})$ known under weak dependence.  This intuition turns out to be
largely incorrect, though a block order $O(n^{1/2})$ may be reasonable (and even optimal) in many cases, owing to  non-linearity  in a long-memory  time series.
While optimal block sizes are more complex under long-range dependence compared to short-range,  we provide a  consistent data-driven rule for block selection, and numerical studies illustrate that the guides for block selection   perform well
in other block-based problems with long-memory time series, such as distribution estimation and  strategies for testing Hermite rank.

\end{abstract}



\end{frontmatter}

\section{Introduction}
\label{sec:1}

Block-based resampling methods  provide useful nonparametric   approximations with statistics from dependent data,
where data blocks help to capture time dependence (cf.~\cite{KL12}).  Considering a stretch from a stationary  series
 $X_1,\ldots,X_n$, a prototypical problem  involves estimating the standard error  of the sample mean $\bar{X}_n=\sum_{t=1}^n X_t/n$.  Subsampling \cite{C86, HJ96, PR94}  and block-bootstrap \cite{K89,LS92}     use sample averages $\bar{X}_{i,\ell}$ computed over length $\ell<n$ data blocks $\{(X_i,\ldots,X_{i+\ell-1})\}_{i=1}^{n-\ell+1}$ within the data $X_1,\ldots,X_n$; in both resampling approaches, the empirical variance of   block averages, say $\widehat{\sigma}_\ell^2$, approximates the block variance  $\sigma_\ell^2\equiv \var(\bar{X}_\ell)$.
If the series $\{X_t\}$ exhibits short-range dependence (SRD) with quickly decaying covariances $r(k)\equiv \cov(X_0,X_k)\rightarrow 0$ as $k\to \infty$ (i.e., $\sum_{k=1}^\infty |r(k)|<\infty$),  then the target variance converges $n \sigma_n^2 =n\var(\bar{X}_n)\rightarrow C>0$ as $n \to \infty$ and $\ell \widehat{\sigma}_\ell^2$ is consistent for $ n\var(\bar{X}_n)$ under mild  conditions ($\ell^{-1}+\ell/n\to \infty$) \cite{PRW99}.
Block-based variance estimators  have further history in time series analysis (cf.~overview in \cite{P03}), including batch means estimation in Markov chain Monte Carlo.     Particularly for SRD, much research has focused on explaining   properties of block-based estimators $\widehat{\sigma}_\ell$  (cf.~\cite{FJ10, K89, L99, L03, PRW99, SS95}).
In turn, these resampling studies have advanced understanding of  best  block sizes (e.g.,  $O(n^{1/3})$) and  implementation under SRD \cite{BK94, HHJ95,   LFL06, NL14, PW04}.

However, in contrast to  SRD, relatively little is known about properties of block-based resampling  estimators  and block sizes under strong  or long-range time dependence (LRD).  For example, recent  tests of Hermite rank \cite{BMG16} as well as other  approximations with block bootstrap and subsampling under LRD  rely on data-blocks \cite{BT17, BTZ16, BW18, JEP16},  creating a need for guides in block selection.  To develop an understanding of data-blocking under LRD, we consider
 the analog problem from SRD of estimating the variance $\var(\bar{X}_n) $ of a sample mean $\bar{X}_n$ through block  resampling
 (cf.~Sec.~\ref{sec:2}-\ref{sec:4}); block selections in this context extend to broader statistics (cf.~Sec.~\ref{sec:5}) and provide  guidance for  distributional approximations with resampling (cf.~Sec.~\ref{sec:6}).
 Because long-memory or long-range dependent (LRD) time series are characterized by
slowly decaying covariances (i.e., $\sum_{k=1}^\infty |r(k)|=+\infty$ diverges),   optimal block sizes in this problem have  intuitively been conjectured as longer $O(n^{1/2})$ than the best block size $O(n^{1/3})$ associated with weak dependence \cite{BW18, HJL98}.  However, this intuition about block selections is misleading.   Under general forms of LRD, the best block selections turn out to depend critically on {\it both} dependence strength (i.e., rate of covariance decay) {\it and} the nature of non-linearity in a time series.
  To illustrate, consider a stationary Gaussian LRD time series $\{Z_t\}$, which we may   associate with  common  models for long-memory \cite{GJ80, MVN68}, and suppose
$\{Z_t\}$ has a covariance decay controlled by a long-memory exponent $\alpha \in (0,1)$ (described more next).
Then, the LRD process $Z_t$ for $\alpha <1/2$ can have an optimal  block length
$O(n^{\alpha})$ while a cousin LRD process
$Z_t+ 0.5 Z_t^2$ has a best block size $O(n^{1/2})$ {\it regardless} of the memory level $\alpha\in (0,1/2)$.
That is,  classes of LRD processes exist  where non-linearity induces a best block order  $O(n^{1/2})$.  Also, as the optimal blocks $O(n^{\alpha})$ for $Z_t$ ($\alpha\in(0,1/2)$) illustrate, when covariance decay slows $\alpha \downarrow 0$ here,
best block sizes for  a resampling variance estimator under LRD do {\it not} generally increase with increasing dependence strength.
While   theory  justifies a block length $O(n^{1/2})$ as optimal in some cases,  the forms of theoretically best block sizes can  generally be complex under LRD and we also establish a provably consistent data-based estimator
of this block size.  Numerical
studies show that the empirical block selection  performs well in variance estimation and provides a guide with good performance in other resampling problems under LRD
 (e.g., distribution estimation for   statistical functionals).

  Section~\ref{sec:2} describes the LRD framework and   variance estimation problem.  We consider  stationary LRD processes $X_t=G(Z_t)$ defined by a transformation $G(\cdot) $ of a LRD Gaussian process $\{Z_t\}$ with a long-memory exponent $\alpha \in (0,1/m)$ (cf.~\cite{T75,T77}); here integer $m \geq 1$
is the so-called Hermite rank of $G(\cdot)$,  which has a well-known impact on the distributional limit of the sample mean $\bar{X}_n$ for such LRD processes (e.g., normal   if $m=1$  \cite{DM79, T79}).
Section~\ref{sec:3} provides the large-sample bias and variance of  block-based resampling   estimators in the sample mean case,
 which are used to determine   MSE-optimal block sizes and a consistent approach to block estimation in Section~\ref{sec:4}.
As complications, best block lengths can depend  on the memory exponent $\alpha$ and
  a higher order rank   beyond Hermite rank $m$ (e.g., 2nd Hermite rank). Two versions of data blocking are also compared, involving fully overlapping (OL) or non-overlapping (NOL) blocks;   while OL blocks are always MSE-better for variance estimation under SRD \cite{L99,L03}, this is not  true under LRD.   
 Section~\ref{sec:5} extends the block resampling   to broader statistical functionals  under LRD (beyond sample means) and includes the block jackknife technique for comparison.    Numerical studies are provided in  Section~\ref{sec:6} to illustrate block selections and resampling
 across problems of variance estimation, distributional approximations, and Hermite rank testing under LRD.
 Section~\ref{sec:7} has concluding remarks, and a supplement \cite{ZLN22} provides proofs and supporting results.

    We end here with some related literature.  Particularly, for Gaussian series $X_t=Z_t$ (or $G(x)=x$ with $m=1$), the computation of  (log)   block-based variance estimators over a series of large block sizes $\ell$ can be a graphical device for estimating the long-memory parameter
    $\alpha$ (using that $\var(\bar{X}_\ell) \approx C_0 \ell^{-\alpha}$ for subsample averages, cf.~Sec.~\ref{sec:2}) \cite{MTT99, TT97}.  Relatedly, \cite{GRS99}  considered block-average regression-type estimators of $\alpha$  in the Gaussian case.  For distribution estimation with LRD linear processes,  \cite{NL05, ZHWW13} studied subsampling, while \cite{KN11} examined
    block bootstrap.  As perhaps the most closely related works, \cite{ADG09,GRS99, KN11}  studied optimal blocks/bandwidths for estimating a sample mean's variance with LRD linear processes (under various assumptions) using data-block averages or related Bartlett-kernel heteroskedasticity
and autocorrelation consistent (HAC) estimators. Those results share connections
   to optimal block sizes here for purely Gaussian series $X_t=Z_t$ (cf.~Sec.~\ref{sec:3.1}), but no empirical estimation rules were considered.   As novelty,
     we account for  LRD data $X_t=G(Z_t)$ from general transformations $G(\cdot)$
   (i.e., the pure Gaussian/linear case is  comparatively simpler), establish consistent block estimation, provide results for more general statistical
   functions, and consider the applications of  block selections in wider contexts under LRD.
   In terms of resampling from LRD transformed Gaussian processes,  \cite{L93} showed the  block bootstrap is valid for approximating the full distribution of the sample mean $\bar{X}_n$ when the Hermite rank is $m =1$,
   while \cite{HJL98} established subsampling as valid for any $m\geq 1$.  (While block bootstrap and subsampling differ in their distributional approximations \cite{PRW99}, these induce a common block-based variance estimator for the sample mean, as described in Sec.~\ref{sec:2}.)
   Recently, much research interest has also focused on further distributional approximations
   with subsampling for   LRD transformed Gaussian processes; see \cite{BT17, BTZ16, BW18, JEP16}.

\section{Preliminaries: LRD processes and block-based resampling estimators}
\label{sec:2}
\subsection{Class of LRD processes} Let $\{Z_t\}$ be a mean zero, unit variance, stationary Gaussian sequence with covariances satisfying
\begin{equation}
\label{eqn:zcov}
\gamma_Z(k)\equiv \E Z_0 Z_k \sim  C_0 k^{-\alpha}
\end{equation}
as $k\to \infty$ for some given $0<\alpha<1$ and constant $C_0>0$.
Examples include fractional Gaussian noise with Hurst parameter $1/2<H<1$ having covariances $\gamma_Z(k)=(|k+1|^{2H} - 2 |k|^{2H} + |k-1|^{2H})/2$
which satisfy (\ref{eqn:zcov}) with $C_0=H(2H-1)$ and $\alpha = 2-2H\in(0,1)$, as well as FARIMA processes with difference parameter $0<d<1/2$ which
satisfy (\ref{eqn:zcov}) with $\alpha= 1-2d \in(0,1)$; see \cite{GJ80, MVN68}.

 Let $G:\mathbb{R}\to \mathbb{R}$ be a real-valued function such that $\E [G(Z_0)]^2<\infty$ holds for a standard normal variable $Z_0$.  In which
 case,  the function $G(Z_0)$ may be expanded
 as
\begin{equation}
\label{eqn:herm}G(Z_0) = \sum_{k=0}^\infty \frac{J_k}{k!}H_k(Z_0)
\end{equation}
  in terms of Hermite polynomials, \[
 H_k(z) \equiv (-1)^k e^{ z^2/2 } \frac{d^k}{dz^k} e^{- z^2/2},\quad k=0,1,2,\ldots,
\]
and corresponding coefficients $J_k \equiv \E [G(Z_0) H_k(Z_0)]$, $k \geq 0$.  The first few Hermite polynomials
are given by $H_0(z)=1$, $H_1(z)=z$, $H_2(z)=z^2-1$, $H_3(z)=z^3-z$, for example, and   $\E H_k(Z_0)=0$ holds for $k\geq 1$.
Let $\mu \equiv \E G(Z_0) = J_0$ denote the mean of $G(Z_0)$ and define the {\it Hermite rank} of $G(\cdot)$  (cf.~\cite{T75}) as
\[
   m  \equiv \min\{ j \geq 1:   J_k \neq 0 \}.
\]
To avoid degeneracy, we assume   $\var[G(Z_0)]>0$  whereby   $m\in [1,\infty)$ is a finite integer.

The target processes of interest are defined as $X_t\equiv G(Z_t)$ with respect to a stationary Gaussian series $ Z_t $ with covariances as in
(\ref{eqn:zcov})
with $0<\alpha<1/m$. Such processes $ X_t $ exhibit strong or long-range dependence
(LRD) as seen by partial sums   $\sum_{k=1}^n |r(k)|$ of covariances
  $r(k)=\cov(X_0,X_k)$ having a slow decay proportional to
\begin{equation}
\label{eqn:1} \sum_{k=1}^n |r(k)| \propto   n^{1-\alpha m} \quad  \mbox{as $n \rightarrow
\infty$, where $  \alpha m\in(0,1)$},
\end{equation}
  depending on the Hermite rank $m$ of the transformation $G(\cdot)$ and memory exponent  $\alpha \in(0,1/m)$ under (\ref{eqn:zcov}).
This represents a common
formulation of LRD, with partial covariance sums diverging as $n \rightarrow
\infty$
 \cite{T75};  see \cite{R95, T02} for further characterizations.

Suppose $X_1,\ldots,X_n$ is an observed time stretch from
the transformed Gaussian series $X_t \equiv G(Z_t)$, having sample mean $\bar{X}_n=n^{-1}\sum_{t=1}^n
X_t$.    Setting $ v_{n,\am}  \equiv
n^{\alpha m}\var(\bar{X}_n)$,   the process structure  (\ref{eqn:zcov})-(\ref{eqn:1})   
entails a so-called long-run variance as
\begin{equation}
\label{xbar.var} \lim_{n \to \infty} v_{n,\am} = v_{\infty,\am} \equiv
  \frac{J_m^2}{m!} \frac{2C_0^m  }{(1-\alpha m)(2-\alpha m)} >0
\end{equation}
(cf.~\cite{ADG09, T75,T77}).   Under LRD,  the variance
$\var(\bar{X}_n)$ of the sample mean decays at a slower rate
$O(n^{-\alpha m})$ as $n \rightarrow \infty$  (i.e., $\alpha m\in(0,1)$) than the typical
$O(n^{-1})$ rate under SRD.    The  limit distribution of $n^{\alpha m/2}(\bar{X}_n-\mu)$ also depends on the Hermite rank $m \geq 1$ \cite{DM79,T77}.

The development   first considers  the variance $ v_{n,\am}  \equiv
n^{\alpha m}\var(\bar{X}_n)$ of the sample mean (or, equivalently here, its limit (\ref{xbar.var}))  as   target of inference under LRD. Resampling results are then extended to broader classes of statistics in Section~\ref{sec:5}.

\subsection{Block-based resampling variance estimators under LRD}
\label{sec2.2}
 A block bootstrap
``re-creates"  the original series $X_1,\ldots,X_n$  by
independently resampling $b \equiv \lfloor n/
\ell\rfloor$  blocks from a collection of length $\ell<n$ data blocks. Resampling from the  overlapping (OL) blocks
$\{(X_i,\ldots,X_{i+\ell-1}):i=1,\ldots,n-\ell+1 \}$ of length $\ell$ within $X_1\ldots,X_n$
yields the moving block bootstrap   \cite{K89, L99, LS92}, while resampling from non-overlapping (NOL) blocks
$\{ (X_{1+\ell(i-1)},\ldots, X_{\ell i}):i=1,\ldots,b \equiv \lfloor n/
\ell\rfloor \}$ gives the
NOL block bootstrap \cite{C86, L03}. Resampled blocks are concatentated to produce a bootstrap series, say $X^*_1,\ldots, X_{\ell b }^*$,
and the distribution of a statistic from the bootstrap series (e.g., $\bar{X}_{\ell b}^* \equiv (\ell b)^{-1}\sum_{i=1}^{\ell b} X_{i}^*$) approximates
 the sampling distribution of an original-data statistic (e.g., $\bar{X}_n$).  Subsampling \cite{PR94, PRW99} is a different approach to approximation 
that computes statistics from  one resampled data block.
Both subsampling and   bootstrap, though,
 estimate a sample mean's variance $v_{n,\am} \equiv n^\am \var(\bar{X}_n)$ with  a common block-based estimator; this is the induced variance of an average under resampling   (e.g., $\var_*(\bar{X}_{\ell b}^*)$), which has a closed  form (cf.~\cite{L93} under LRD), resembling a batch means estimator \cite{FJ10}.
Based on $X_1,\ldots,X_n$, the OL block-based variance estimator of $v_{n,\am} \equiv n^\am \var(\bar{X}_n)$ is given by
\begin{equation}
\label{boot.var}
\widehat{V}_{\ell,\am,\OL} =
\frac{1}{n-\ell+1}
\sum_{i=1}^{n-\ell+1}\ell^\am(\bar{X}_{i,\ell}-\widehat{\mu}_{n,\OL})^2,
\quad \widehat{\mu}_{n,\OL}=\frac{1}{n-\ell+1}\sum_{i=1}^{n-\ell+1}\bar{X}_{i,\ell},
\end{equation}
where above $\bar{X}_{i,\ell}= \sum_{j=i}^{i+\ell-1} X_j/\ell$ is the sample average of the $i$th data block $(X_i,\ldots,X_{i+\ell-1})$,   $i=1,\ldots,n-\ell+1$.  Essentially, block versions $\{\ell^{\am/2} (\bar{X}_{i,\ell} -\bar{X}_n)\}_{i=1}^{n-\ell+1}$ of  the quantity $n^{\am/2}
(\bar{X}_n-\mu)$ give a sample variance $\widehat{V}_{\ell,\am,\OL}$ that estimates
$v_{\ell,\am} \equiv \ell^\am \var(\bar{X}_\ell) \approx v_{n,\am} \equiv n^\am \var(\bar{X}_n) $ for  sufficiently large $\ell,n$ by (\ref{xbar.var}).
The  NOL block-based  variance estimator is
 defined as
\[
\widehat{V}_{\ell,\am,\NOL}  =
\frac{1}{b}\sum_{i=1}^{b} \ell^\am(\bar{X}_{1+(i-1)\ell,\ell} -
\widehat{\mu}_{n,\NOL})^2, \quad \widehat{\mu}_{n,\NOL}=\sum_{i=1}^{b}
\bar{X}_{1+(i-1)\ell,\ell}/b,
\]
using NOL    averages $\bar{X}_{1+(i-1)\ell,\ell}$, $i=1,\ldots,b  \equiv \lfloor n/ \ell\rfloor$, where $\widehat{\mu}_{n,\NOL} = \bar{X}_n$ when $n=b\ell$.

Under SRD,   variance estimators are standardly defined by letting $\am =1$ above (e.g., in $\widehat{V}_{\ell,\am,\OL}$ from (\ref{boot.var})).   Likewise, under LRD, both the target variance
$v_{n,\am} \equiv n^\am \var(\bar{X}_n)$  and block-based estimators
$\widehat{V}_{\ell,\am,\OL}$ or $\widehat{V}_{\ell,\am,\NOL}$ are scaled   to be comparable,  which involves the long-memory exponent $\am\in(0,1)$.
In practice, $\am\in(0,1)$ is usually unknown.  To develop  block-based estimators under LRD, we first consider
  $\am \in (0,1)$ as given. Ultimately, an estimate $\widehat{\am}_n$ of $\am$ can be substituted into
    $\widehat{V}_{\ell,\am,\OL}$
or $\widehat{V}_{\ell,\am,\NOL}$ which, under mild conditions, does not change conclusions about consistency or  best estimation rates  (cf.~Sec.~\ref{sec:4.2}).

\section{Properties for block-based resampling estimators under LRD}

\label{sec:3}

Large-sample results for the block-based variance estimators   require some  extended   notions of the Hermite rank  of $G(\cdot)$ in defining $X_t \equiv G(Z_t) = \sum_{k=0}^\infty J_k/k!
\cdot H_k(Z_t)$, for $J_k = \E G(Z_0) H_k(Z_0)$, $k \geq 1$. Recalling
 $m  \equiv \min\{ k >0:  J_k \neq 0   \}$ as the usual   Hermite rank   of $G(\cdot)$,  define the {\it 2nd Hermite rank} of $G(\cdot)$ by the index
\[
 m_2 \equiv \min\{ k > m:  J_k \neq 0   \}
\]
of the next highest non-zero coefficient in the Hermite expansion
 (\ref{eqn:herm}) of $G(\cdot)$.  In other words, $m_2$ is the Hermite rank of
$X_t -\mu - J_{m} H_{m}(Z_t)/m!$  upon removing the mean and 1st Hermite rank term  from $X_t=G(Z_t)$.  If the set
$\{ k > m:  J_k \neq 0   \}$ is empty, we define
$m_2 = \infty$.   We also define the {\it Hermite pair-rank} of a function $G(\cdot)$ by
\[
m_p \equiv \inf\{  k \geq m:  J_k J_{k+1} \neq 0   \};
\]
when the above set is empty, we define $m_p =\infty$.
The Hermite
pair-rank $m_p$ identifies the index of the first consecutive pair of non-zero terms $(J_k,J_{k+1})$  in the expansion $X_t=G(Z_t)=\mu + \sum_{k=1}^\infty J_k/k! H_k(Z_t) $.  For   non-degenerate series
$X_t=G(Z_t)$, the Hermite rank $m$ is always finite, but  the 2nd rank $m_2$ and   pair-rank $m_p$ may not be (and $m_2=\infty$ implies $m_p=\infty$).  For example, both   series $G(X_t) =  H_1(Z_t)$ and $G(X_t) =  H_1(Z_t)+H_{3}(Z_t)$   have Hermite rank $m=1$, pair-rank
$m_p=\infty$, and 2nd  ranks $m_2=\infty$ and $3$, respectively; the series  $G(X_t) = \ H_1(Z_t)+H_3(Z_t)+H_4(Z_t)$ and
$G(X_t) =  H_3(Z_t)+H_4(Z_t)$   have pair-rank $m_p=3$  with Hermite ranks $m=1$ and $3$, and 2nd  ranks $m_2=3$ and $4$, respectively.

In what follows, due to combined effects of dependence and non-linearity in a LRD time series $X_t = G(Z_t)$,  the Hermite pair-rank $m_p \in[m,\infty]$
of $G$ plays a role in the asymptotic variance   of   resampling estimators (Sec.~\ref{sec:3.2}), while
the 2nd Hermite rank $m_2\in [m+1,\infty]$ impacts the bias of
resampling estimators (Sec.~\ref{sec:3.1}).

\subsection{Large-sample bias  properties}
\label{sec:3.1}
Bias expansions for the block resampling estimators require a
more detailed form of the LRD covariances than (\ref{eqn:zcov}) and we
suppose that
\begin{equation}
\label{eqn:r} \gamma_Z(k)\equiv \cov(Z_0,Z_k) =  C_0 k^{-\alpha}(1 + k^{-\tau} L(k)), \quad k >0,
\end{equation}
holds for some $\alpha \in(0,1/m)$ and $C_0>0$ (again $\gamma_Z(0)=1$) with some  $\tau \in  (1-\alpha m,\infty) $ and  slowly varying function
$L: \mathbb{R}^+ \to \mathbb{R}^+$  that  satisfies $\mathrm{lim}_{x\to \infty}   L(ax)/L(x)  = 1$ for any $a>0$.   For Gaussian FARIMA (i.e., $\alpha=1-2d \in (0,1)$) and Fractional Gaussian noise (i.e., $\alpha=2-2H \in(0,1)$) processes $\{Z_t\}$, one
may verify that (\ref{eqn:r}) holds with $\tau=1$ for any $\alpha\in(0,1/m)$ and $ m \geq 1$.
 The statement of bias in Theorem~\ref{theorem1} additionally requires process constants $\mathcal{B}_0(m),\mathcal{B}_1(m_2)$ that depend on the  1st $m$ and 2nd $m_2$ Hermite ranks and covariances  in (\ref{eqn:r}).  These are given as $\mathcal{B}_1(m_2)\equiv 2  \sum_{j=m_2}^\infty  (J^2_j/j!)
      \sum_{k=1}^\infty   [\gamma_Z(k)]^j$ when $\alpha m_2> 1$; $\mathcal{B}_1(m_2)\equiv     2  C_0^{m_2} J_{m_2}^2 /m_2!$
  when $\alpha m_2 =  1$;   $\mathcal{B}_1(m_2)\equiv 2 J_{m_2}^2 C_0^{m_2}/[{m_2}! (1-m_2\alpha)(2-m_2 \alpha) ] = v_{\infty,\alpha m_2}$ from (\ref{xbar.var}) when $\alpha m_2 <  1$; and
\begin{equation}
\label{eqn:B1}
\mathcal{B}_0(m) \equiv    2  C_0^m \frac{J_m^2}{m!}\left\{\mathcal{I}_\am + \sum_{k=1}^\infty  k^{-\am} \sum_{j=1}^m {m \choose j} [L(k) k^{-\tau}]^j  \right\}   +
\sum_{k=m}^\infty  \frac{J_k^2}{k!},\end{equation}
with Euler's generalized constant
$\mathcal{I}_\am \equiv  \lim_{k\to \infty} ( \sum_{j=1}^{k} j^{-\alpha m} - \int_{0}^k x^{-\alpha m} dx ) \in (-\infty, 0)$.

\begin{theorem}
\label{theorem1} {\it Suppose  $X_t \equiv G(Z_t)$  where the stationary Gaussian process $\{Z_t\}$ satisfies
(\ref{eqn:r}) with memory exponent $\alpha \in (0,1/m)$ and where $G(\cdot)$ has Hermite rank $m$ and 2nd  Hermite rank $m_2$ (noting $m_2 > m$ and possibly $m_2=\infty$}).
  Let $\widehat{V}_{\ell,\am}$ denote either
$\widehat{V}_{\ell,\am,\OL}$ or
$\widehat{V}_{\ell,\am,\NOL}$ as block resampling estimators of $v_{n,\am} = n^{\am} \var(\bar{X}_n)$ based on $X_1,\ldots,X_n$.  If  $\ell^{-1}+
\ell/n \rightarrow 0$ as $n \rightarrow \infty$, then    the bias of $\widehat{V}_{\ell,\am}$
is given by
\begin{eqnarray*}
  \E\widehat{V}_{\ell,\am} - v_{n, \am} &= &\mathcal{B}_{0}(m)\left(\frac{\ell^{\am}}{\ell}\right)\big(1+o(1)\big) -  v_{\infty,\am}
  \left(\frac{\ell}{n}\right)^{\am}  \big(1+o(1)\big)\\ &&\quad +I(m_2<\infty)\mathcal{B}_{1}(m_2)\left(\frac{\ell^{\am}}{\ell^{\min\{1, \alpha
  m_2\}}}\right)[\log \ell]^{I(\alpha m_2=1)}\big(1+o(1)\big),
\end{eqnarray*}
where $I(\cdot)$ denotes the indicator function, the constant $v_{\infty,\am} \equiv  2 J_m^2 C_0^m/[m! (1-\alpha m)(2-\alpha m) ]$  is from (\ref{xbar.var}), and constants $\mathcal{B}_{0}(m),\mathcal{B}_1(m_2)$ are from (\ref{eqn:B1}).
\end{theorem}
\noindent {\sc Remark 1:} If we switch the target of variance estimation  from $v_{n,\am} = n^{\am} \var(\bar{X}_n)$ to the limit variance
$v_{\infty,\am} \equiv \lim_{n\to \infty} v_{n,\am}$ from (\ref{xbar.var}), this does not change the bias expansion in Theorem~\ref{theorem1}
or  results in Section~\ref{sec:4} on best block sizes for minimizing MSE.   \\

To better understand the  bias   of a block-based estimator  under LRD, we may consider the   case of a purely   Gaussian LRD  series $X_t = Z_t$ (i.e., no transformation), corresponding to $G(x)=x$, $m=1$ and $m_2 = \infty$.   The
bias
 then    simplifies under Theorem~\ref{theorem1} as
\begin{equation}
\label{eqn:Lin}
  \E\widehat{V}_{\ell,\am} - v_{n, \am} =   \mathcal{B}_0(1)  \left(\frac{\ell^\alpha}{\ell}\right) \big(1+o(1)\big) -  v_{\infty,\am}
  \left(\frac{\ell}{n}\right)^{\alpha}  \big(1+o(1)\big),
\end{equation}
  depending only on the memory exponent $\alpha$ of the process $ Z_t $. This bias form can also hold when $X_t$ is
   LRD and linear  \cite{GRS99,KN11}.  However, for a transformed LRD series $X_t=G(Z_t)$,  the function $G$   and the underlying   exponent $\alpha$
  impact   the bias of the block-based   estimator in intricate ways.
  The order of a main bias term in Theorem~\ref{theorem1} is generally summarized as
\begin{equation}
\label{eqn:biasO}
O\left(\frac{ [\log \ell]^{I(\alpha m_2 =1)}}{\ell^{\min\{1,\alpha m_2\} - \am}}\right),
\end{equation}
 which depends on how the 2nd  Hermite rank $m_2>m$ of the transformed series $X_t \equiv G(Z_t)$, as a type of non-linearity measure,
 relates to the long-memory exponent $\alpha \in (0,1/m)$.  Small
values of $m_2$ satisfying $1/\alpha > m_2$  induce  the worst bias rates $O(\ell^{- (m_2-m)\alpha})$   compared to the
 best  possible bias $O(\ell^{-(1-\am)})$  occurring, for example,  when $m_2=\infty$ (or   no terms in the Hermite
expansion of $G(\cdot)$ beyond the 1st rank $m$). In fact, the largest bias rates arise whenever  2nd Hermite rank terms $J_{m_2} H_{m_2}(Z_t)/m_2!$
exist in the expansion of $X_t \equiv G(Z_t)$ (i.e., $J_{m_2} \neq 0$) and  exhibit long-memory
under $\alpha m_2 <1$. (cf.~Sec.~A.1 in \cite{ZLN22}).    For comparison,   block-based    estimators in the SRD case  \cite{L99} exhibit a  smaller bias  $O(1/\ell)$
  than   the best  possible bias  in (\ref{eqn:biasO}) under LRD.

\subsection{Large-sample variance  properties}
\label{sec:3.2}
To   establish the variance of the block resampling   estimators under LRD, we require an additional moment condition regarding the transformed series $X_{t}=G(Z_t)$.   For second   moments, a simple characterization exists that
$\E X_t^2 = \E [G(Z_t)]^2<\infty$ is finite if and only if $\sum_{k=0}^\infty J_k^2/k! <\infty$.
 For higher order moments, however,  more elaborate conditions are required to guarantee $\E X_t^4 = \E [G(Z_t)]^4<\infty$ and perform   expansions of $\E X_{t_1} X_{t_2} X_{t_3} X_{t_4}$.  We shall use a condition ``$G\in \overline{\mathcal{G}}_4(1)$" from \cite{T77}.
(More generally, Definition~3.2 of \cite{T77} prescribes a condition $G\in \overline{\mathcal{G}}_4(\epsilon)$, with $\epsilon \in (0,1]$,  for moment expansions,
which could   be applied to derive Theorem~\ref{theorem2} next.  We use $\epsilon=1$   for simplicity, where a sufficient condition for $G\in \overline{\mathcal{G}}_4(1)$  is
   $\sum_{k=0}^\infty  3^{k/2} |J_k| /\sqrt{k!}<\infty$, holding for any polynomial $G$, cf.~\cite{T77}).  See the supplement~\cite{ZLN22}  for
more technical  details.

To state the large-sample variance properties of block-based estimators $\widehat{V}_{\ell,\am,\OL}$ or
$\widehat{V}_{\ell,\am,\NOL}$ in Theorem~\ref{theorem2}, we also require some proportionality constants.
 As a function of the Hermite rank, when $m \geq 2$ and $ \alpha m<1  $,   define a positive scalar
 \begin{equation*}
    \phi_{\alpha,m}\equiv  \frac{2}{(1-2\alpha  )(1-\alpha)} \left( \frac{2 J_m^2C_0^m}{ (m-1)!}\frac{1}{[1-(m-1)\alpha][2-(m-1)\alpha)]} \right)^2.
    \end{equation*}
In the case of a Hermite rank $m=1$, define another    positive proportionality constant, as a function of $\alpha\in(0,1)$ and the type  of resampling blocks (OL/NOL)), as
   \[
 a_{\alpha}
 \equiv \frac{8 J_1^4 C_0^2  }{(1-\alpha)^2(2-\alpha)^2} \times \left\{\begin{array}{cll} 1 +
\frac{(2-\alpha)^2(2\alpha^2+3\alpha-1)}{4(1-2\alpha)(3-2\alpha)}-
\frac{ \Gamma^2(3-\alpha)}{\Gamma(4-2\alpha)} &\,\, &
 \mbox{if
$0<\alpha<1/2$, OL or
NOL}\\[.6cm]
  9/32&&\mbox{if $\alpha=1/2$, OL or NOL}\\[.3cm]
  \sum_{x = -\infty}^\infty g^2_\alpha(x) & &\mbox{if $1/2<\alpha < 1$, NOL}\\[.7cm]
 \int_{-\infty}^\infty g^2_\alpha(x) dx  & &\mbox{if $1/2<\alpha < 1$, OL,}\\
\end{array}
 \right.
\]
where $\Gamma(\cdot)$ denotes the gamma function and
$g_\alpha(x)\equiv(|x+1|^{2-\alpha}-2|x|^{2-\alpha} +
|x-1|^{2-\alpha})/2$, $x \in \mathbb{R}$. In the definition of
$a_{\alpha} $,
 $g_\alpha^2(x)$ is summable/integrable when
$\alpha\in(1/2,1)$  using  $g_\alpha(x) \sim (2-\alpha)(1-\alpha)
x^{-\alpha}/2$  as $x \rightarrow \infty$.  Finally, as a function of any   Hermite pair-rank    $m_p\in [1,\infty]$ and $\alpha \in
(0,1)$, we define
a constant as
\begin{equation*}
\lambda_{\alpha,m_p} \equiv\frac{8 C_0}{(1-\alpha)(2-\alpha)}  \times \left\{\begin{array}{cl}
 (\frac{2 C_0^{m_p}J_{m_p}J_{m_p+1}}{m_p!})^2 [(1-\alpha m_p)(2-\alpha m_p)]^{-2I(\alpha m_p <1)}    &\mbox{if $\alpha m_p \leq 1$}\\[.3cm]
 \left(\sum_{k=m_p}^\infty  \sum_{j=-\infty}^\infty[\gamma_Z(j)]^k  J_k J_{k+1}/k!\right)^2  &\mbox{if $\alpha m_p > 1$},
\end{array}
 \right.
\end{equation*}
with Gaussian covariances $\gamma_Z(\cdot)$ and $C_0>0$ from (\ref{eqn:zcov})  and an indicator $I(\cdot)$ function.

With constants $\lambda_{\alpha,m_p},\phi_{\alpha,m}, a_{\alpha} > 0$ as above, we may next state Theorem~\ref{theorem2}.

\begin{theorem}
\label{theorem2}
{\it Suppose  $ X_t \equiv G(Z_t) $  where the stationary Gaussian process $\{Z_t\}$ satisfies
(\ref{eqn:zcov}) with  $C_0>0$ and memory exponent $\alpha \in (0,1/m)$ and where  $G\in \overline{\mathcal{G}}_4(1)$  has Hermite rank $m \geq
1$    and Hermite pair-rank $m_p$ (note $m_p\geq m$ and possibly $m_p=\infty$}). Let $\widehat{V}_{\ell,\am}$ denote either
$\widehat{V}_{\ell,\am,\OL}$ or
$\widehat{V}_{\ell,\am,\NOL}$ as block resampling estimators of $v_{n,\am} = n^{\am} \var(\bar{X}_n)$ based on $X_1,\ldots,X_n$.  If  $\ell^{-1}+
\ell/n \rightarrow 0$ as $n \rightarrow \infty$, then    the variance of $\widehat{V}_{\ell,\am}$ is given   by
\[
\var(\widehat{V}_{\ell,\am}) = \left\{  \begin{array}{cll}
\displaystyle{ \phi_{\alpha,m} \left(\frac{\ell}{n}\right)^{2\alpha} \big(1+o(1)\big) + r_{n,\alpha,m,m_p} } &\,\,& \mbox{if $m \geq 2$}\\[.4cm]
\displaystyle{  a_\alpha \left(\frac{\ell}{n}\right)^{\min\{1,2\alpha\}} [\log n ]^{I(\alpha =1/2)} \big(1+o(1)\big) + r_{n,\alpha,m,m_p}} &\,\,& \mbox{if $m=1$},
\end{array}
\right. \]
where $I(\cdot)$  denotes an indicator function  and
\[
r_{n,\alpha,m,m_p}
\equiv I(m_p<\infty)\frac{\lambda_{\alpha,m_p}}{n^\alpha}
\left(  \frac{\ell^{\am}}{\ell^{\min\{1,\alpha m_p\}}}  [\log \ell]^{I( \alpha m_p=1)}\right)^2
\big(1+o(1)\big).
\]

\end{theorem}
\noindent {\sc Remark 2:} Above $r_{n,\alpha,m,m_p}$ represents a second variance contribution, which depends on the Hermite pair-rank $m_p$ and is non-increasing in block length $\ell$ (by  $1/\alpha<m\leq m_p$). The value of $r_{n,\alpha,m,m_p}$  is  zero when $m_p=\infty$ and is largest $O(n^{-\alpha})$ when the pair-rank  assumes its smallest possible value $ m_p=m$. For example, the series $X_t=  H_m(Z_t)+H_{m+1}(Z_{t})$ and  $X_t=  H_m(Z_t)+H_{m+2}(Z_{t})$  have pair-ranks $m_p=1$ and $\infty$, respectively, inducing different $r_{n,\alpha,m,m_p}$ terms.    While  $r_{n,\alpha,m,m_p}$ can  dominate the variance expression of Theorem~\ref{theorem2} for some block $\ell$ sizes,  the contribution of $r_{n,\alpha,m,m_p}$  emerges as asymptotically negligible at an optimally selected block size $\ell_{opt}$ (cf.~Sec.~\ref{sec:4.1}).\\

By Theorem~\ref{theorem2}, the variance  of a resampling estimator  $\widehat{V}_{\ell,\am}$  depends on the block size $\ell$  through a decay rate  $O((\ell/n)^{2
\alpha})$ that, surprisingly, does not involve the exact value of the rank Hermite $m$.  The reason is that, when $ m \geq 2$,    fourth order cumulants
of the transformed process $X_t = G(Z_t)$ determine this variance (cf.~\cite{ZLN22}).
Also, any differences in block type ($\widehat{V}_{\ell,\am,\OL}$ vs
$\widehat{V}_{\ell,\am,\NOL}$) only emerge in a proportionality  constant $a_{\alpha}$
when $m=1$ and $\alpha \in (1/2,1)$; otherwise,  $a_{\alpha}$ does not change with block type.
 Consequently, for  processes $X_t=G(Z_t)$ with  strong LRD ($\alpha<1/2$, $m\geq 1$), there is no large-sample advantage to  OL   blocks for  variance estimation. 
 In contrast,  under SRD,
OL blocks  reduce the  variance of a resampling variance estimator by a multiple of $2/3$ compared to  NOL blocks  \cite{K89,L99,L03}, because the non-overlap between two OL blocks (e.g., $X_1,\ldots,X_{\ell}$ and $X_{1+i},\ldots,X_{\ell+i}$, $i<\ell$)   acts 
 roughly uncorrelated. This fails under strong LRD  where
OL/NOL blocks have the same variance/bias/MSE properties here. Section~\ref{sec:6} provides  numerical examples. As under SRD, however, OL blocks remain generally preferable (i.e., smaller $a_\alpha$ for weak LRD. $\alpha>1/2$).


\section{Best Block Selections and Empirical Estimation}
\label{sec:4}
\subsection{Optimal Block Size and MSE}
\label{sec:4.1}
 Based on the large-sample bias and variance expressions in Section~\ref{sec:3},
an explicit form for the optimal block size $\ell_{opt}\equiv \ell_{opt,n}$  can be determined for minimizing  the asymptotic mean squared error
\begin{equation}
\label{eqn:MSE1}
\mathrm{MSE}_{n}(\ell)  \equiv \E( \widehat{V}_{\ell,\am} - v_{n,\am})^2
 \end{equation}
 of a block-based resampling estimator  $\widehat{V}_{\ell,\am}$  of $v_{n,\am}\equiv n^{\am} \var(\bar{X}_n)$   under LRD.
 \begin{cor}
 \label{cor1}
Under   Theorems~\ref{theorem1}-\ref{theorem2} assumptions, the optimal block size for a resampling estimator $\widehat{V}_{\ell,\am,\OL}$ or $\widehat{V}_{\ell,\am,\NOL}$    is given by (as $n\to\infty$)
\begin{equation*}
 \ell_{opt,n} =K_{\alpha,m,m_2}\times
\begin{cases}
 \displaystyle{n^{ \frac{\alpha}{\alpha(1 - m) + \min\{1,\alpha m_2\}}}(\log n)^{I(\alpha m_2=1)}} & \text{if } 0 < \alpha < 0.5, \; m\geq 1 \\
  (\frac{n}{\log n})^{0.5} (\log n)^{I(\alpha m_2=1)} & \text{if } \alpha = 0.5, \; m=1 \\
  n^{\frac{1}{3 - 2\alpha}}& \text{if } 0.5 < \alpha < 1, \; m=1, \\
\end{cases}
\end{equation*}
for a constant $K_{\alpha,m,m_2}>0$,  changing by block type OL/NOL only when $m=1,\alpha\in(1/2,1)$.
\end{cor}
The Appendix provides values for $K_{\alpha,m,m_2}>0$.
For LRD processes $X_t=G(Z_t)$, best block lengths $\ell_{opt}$  depend intricately on the transformation $G(\cdot)$ (through   ranks $m,m_2$)    and the memory parameter $\alpha  < 1/m$  of the Gaussian process $ Z_t $.  Optimal blocks  {\it increase} in length whenever the strength of long-memory {\it decreases} (i.e., $\alpha$ increases); as $\alpha$ moves closer to $1/m$,  the order of $\ell_{opt}$   moves closer to $O(n)$. This is a counterintuitive aspect of LRD in resampling.  With
  variance estimation under SRD \cite{L99,L03}, best block size has a known order $ C n^{1/3}$ where the proportionality constant $C>0$  increases  with dependence.

The 2nd Hermite rank $m_2$ of $G(\cdot)$ can particularly impact   $\ell_{opt}$.  Whenever $\alpha<1/m_2$,  the optimal block order $\ell_{opt}\propto n^{1/(m_2-m+1)}$ does  not change.  As a consequence in this case, if an immediate second term $H_{m+1}(Z_0)$ appears in the Hermite expansion (\ref{eqn:herm}) of $X_0=G(Z_0)$, so that the 2nd  rank  is $m_2=m+1$, then
  the optimal block size   becomes $\ell_{opt}\propto  n^{1/2}$.
This suggests that a guess $\ell_{opt} = O(n^{1/2})$ often found in the literature for block resampling under LRD  can  be reasonable, though not by the   intuition
that slow covariance decay under LRD   implies larger blocks compared to those $O(n^{1/3})$ for SRD. Rather, for
transformations $G(\cdot)$ where $m_2=m+1$ may hold naturally, the choice $\ell_{opt}\propto n^{1/2}$ is optimal with sufficiently strong $\alpha<1/(m+1)$ dependence, regardless of the exact Hermite rank $m$.

For completeness, we note that $\mathrm{MSE}_{n}(\ell)  \equiv \E( \widehat{V}_{\ell,\am} - v_{n,\am})^2$  has an optimized   order as
\begin{equation}
\label{eqn:MSE2}
 \mathrm{MSE}_{n}(\ell_{opt,n})   \propto
\begin{cases}
 n^{ \frac{-2\alpha(\min\{1,\alpha m_2\} - \am  )}{\alpha(1 - m) + \min\{1,\alpha m_2\}}} (\log n)^{2\alpha I(\alpha m_2=1)} & \text{if } 0 < \alpha < 0.5, \; m\geq 1 \\
 (\frac{n}{\log{n}})^{-0.5}  (\log n)^{I(\alpha m_2=1)}   & \text{if } \alpha = 0.5, \; m=1 \\
 n^{- \frac{2(1 - \alpha)}{3 - 2 \alpha}}  & \text{if } 0.5 < \alpha < 1, \; m=1. \\
\end{cases}
\end{equation}
at the optimal  block $\ell_{opt}\equiv \ell_{opt,n}$, which also depends on   $m,m_2$ and $\alpha$    under LRD.

Section~\ref{sec:4.2} next shows that estimation of the long-memory exponent does not affect the large-sample results
and block considerations for the resampling variance estimators.  Section~\ref{sec:4.3} then provides
a consistent data-driven method for estimating the  block size $\ell_{opt}$.


\subsection{Empirical considerations for long-memory exponent}
\label{sec:4.2}
We  have assumed  the memory exponent $\am \in (0,1)$ of the LRD process $X_t=G(Z_t)$ is known in   block-based
resampling estimators  (\ref{boot.var}).  If an appropriate estimator $\widehat{\am}_n$ of  $\am$ is
instead substituted,  then the resulting
estimators will possess similar consistency rates under mild conditions.

We let $ \widehat{V}_\ell$
 generically denote a  block-based estimator
$\widehat{V}_{\ell,\am}$ of $v_{n,\am}=n^{\am}\var(\bar{X}_n)$ (e.g., $\widehat{V}_{\ell,\am,\OL}$
or  $\widehat{V}_{\ell,\am,\NOL}$)  found by replacing
$\am$ with an estimator $\widehat{\am}_n$ based on $X_1,\ldots,X_n$.

\begin{cor}
\label{cor2}  Suppose Theorem~\ref{theorem1}-\ref{theorem2} assumptions.
 As   $n\to \infty$,\\[.1cm]
(i) if $| \widehat{\am}_n-\am| \log n \stackrel{p}{\rightarrow}0$, then  $\widehat{V}_{\ell}$ is consistent for
$v_{n,\am}$,  i.e., $|\widehat{V}_{\ell} - v_{n,\am}| \stackrel{p}{\rightarrow}0$.\\
(ii)  if $| \widehat{\am}_n-\am| \log n =O_p(n^{-1/4})$, then the convergence rate of $|\widehat{V}_{\ell}- v_{n,\am}|$
in probability   matches that of $|\widehat{V}_{\ell,\am} - v_{n,\am}|$ from  Theorems~\ref{theorem1}-\ref{theorem2}.
\end{cor}
Several potential estimators of $\am$ satisfy  Corollary~\ref{cor2} conditions, such as log-periodogram or local Whittle estimation
\cite{R95}.  These, for example, can exhibit sufficiently fast convergence in probability, e.g. $O_p(n^{-2/5})$ (cf.~\cite{AS04,HDB98}).
For simplicity, we use local Whittle estimation with bandwidth $\lfloor n^{0.7} \rfloor$  (cf.~\cite{AS04}) in the following.

\subsection{Data-driven block estimation}
\label{sec:4.3}
The   block results from Section~\ref{sec:4.1} suggest that
  data-driven choices of block size under LRD have no simple analogs to block-resampling in the SRD case.
  For variance estimation under SRD,
  several approaches exist
    for estimating the best   block size $\ell_{opt} $ through    plug-in estimation \cite{BK94,  LFL06, PW04} or
empirical MSE-minimization (\cite{HHJ95}-method).  By exploiting   the known block order $\ell_{opt} \approx C n^{1/3}$ under SRD, these methods
  target the proportionality constant $C>0$.     In contrast,
optimal blocks under LRD have a form
$\ell_{opt } \approx   K n^{a} (\log n )^{\tilde{I}}$ from Corollary~\ref{cor1}, where $K\equiv K_{\alpha, m, m_2}>0$ and $a\equiv a_{\alpha, m, m_2} \in(0,1)$ are complicated terms based on $\alpha, m, m_2$, while   $\tilde{I} \equiv -0.5  I(\alpha=0.5, m  = 1)+ I(\alpha m_2 = 1)$ involves indicator functions. Because the  order $n^{a} (\log n )^{\tilde{I}}$
 is   unknown  in practice,   previous strategies to block estimation are not directly applicable in the LRD setting.
   Plug-in estimation  seems particularly intractable under LRD;   general plug-in approaches  under SRD \cite{LFL06}
 require   known orders for
    bias/variance  in   estimation, but these are also unknown under LRD
   (Theorems~\ref{theorem1}-\ref{theorem2}).
 Consequently, we consider a modified method for estimating block size $\ell_{opt}\equiv \ell_{opt,n}$   that involves {\it two} rounds of empirical MSE-minimization (\cite{HHJ95}-method)

To adapt the \cite{HHJ95}-method for LRD, we take a collection of subsamples $(X_i,\ldots,X_{i+h-1})$ of length $h<n$, $i=1,\ldots,n-h+1$.  Based on $X_1\ldots,X_n$, let  $\widehat{\am}_n$ denote an estimator of $\am$ (for use in all estimators to follow)   and let $\widehat{V}_{\tilde{\ell}}$ denote a resampling variance estimator   (replacing   $\am$ with $\widehat{\am}_n$) based on a pilot block size $\tilde{\ell}$ (e.g., $\tilde{\ell} \propto n^{1/2}$).
Similarly, let $\widehat{V}_{ \ell}^{(i,h)}  $ denote a resampling variance estimator computed on the subsample $(X_i,\ldots,X_{i+h-1})$ using a block length $\ell<h$, $i=1,\ldots,n-h+1$.  We then define an initial block-length estimator $\widehat{\ell}_{opt, h}$ as the
minimizer of the empirical MSE
\[
\widehat{\mathrm{MSE}}_{h,n}(\ell)   \equiv \frac{1}{n-h+1}\sum_{i=1}^{n-h+1} \left( \widehat{V}_{ \ell}^{(i,h)} -
\widehat{V}_{\tilde{\ell}} \right)^2,\quad   1 \leq \ell < h.
\]
Here $\widehat{\mathrm{MSE}}_{h,n}(\ell)$ estimates  $\mathrm{MSE}_h(\ell)$ in (\ref{eqn:MSE1}), or the MSE of a resampling estimator based on a sample of size $h$ and block size $\ell$, while $\widehat{\ell}_{opt, h}$ then estimates the minimizer of $\mathrm{MSE}_h(\ell)$ or the optimal block ${\ell}_{opt, h}$ from Corollary~\ref{cor1} using ``$h$" in place of ``$n$" there.
 Above the pilot estimator $\widehat{V}_{\tilde{\ell}}$
plays the role of a target variance to mimic the MSE formulation (\ref{eqn:MSE1}).

Theorem~\ref{cor3} establishes important conditions on the subsample size $h$ and pilot block $\tilde{\ell}$ for consistent estimation under LRD.
For the   transformed series $X_t=G(Z_t)$, the result involves a general 8th order moment condition (i.e., $G\in \overline{\mathcal{G}}_8(1)$ under Definition~3.2 of \cite{T77}) analogous  to the 4th order moment condition described in Section~\ref{sec:3.2}.

\begin{theorem}
\label{cor3} Along with Theorem~\ref{theorem1}-\ref{theorem2} assumptions, suppose $| \widehat{\am}_n-\am| \log n =O_p(n^{-1/4})$ (as in Corollary~\ref{cor2}) and that $G\in \overline{\mathcal{G}}_8(1)$.  Suppose also that $h\to \infty$ and $\tilde{\ell}\to \infty$   with $h/\tilde{\ell} + \tilde{\ell} h/n = O(1)$. Then, as $n\to \infty$, the empirical MSE, $\widehat{\mathrm{MSE}}_{h,n}(\ell)$,  has a sequence of minimizers $\widehat{\ell}_{opt, h}$ such that
\[
\frac{\widehat{\ell}_{opt, h}}{ \ell_{ opt,h }}  \overset{p}{\to} 1 \quad \mbox{and}\quad \frac{\widehat{\mathrm{MSE}}_{h,n}(\widehat{\ell}_{opt, h})}{{\mathrm{MSE}}_{h}( \ell_{opt, h})} \overset{p}{\to} 1,
\]
where   $ \ell_{ opt,h } $ has Corollary~\ref{cor1} form  with $\mathrm{MSE}_{h}( \ell_{opt, h})$ as in  (\ref{eqn:MSE2}) (with ``$h$" replacing ``$n$").
\end{theorem}
 Theorem~\ref{cor3} does not address  estimation of the best block size $\ell_{opt,n}$ for a length $n$ time series,   but rather the optimal block $\ell_{opt, h}$ for a smaller length $h<n$   series.  Nevertheless,
 the result establishes a  non-trivial first step that, under LRD, some block sizes  can be validly estimated  through empirical MSE (\cite{HHJ95}-method) provided that
the subsample size $h$ and pilot block $\tilde{\ell}$ are appropriately chosen.  In particular, the condition $h/\tilde{\ell} + \tilde{\ell} h/n = O(1)$
cannot be reduced (related to pilot estimation $\widehat{V}_{\tilde{\ell}}$) and entails that the largest subsample length possible is $h=O(n^{1/2})$ within the empirical MSE approach under LRD.

With this in mind and because the order of $\ell_{opt,n}$ is unknown,
we  use empirical MSE device twice, based on two  subsample lengths $h \equiv h_1=\lfloor C_1 n^{1/r}\rfloor$ and $h_2 =\lfloor C_2 n^{\theta/r} \rfloor$. Here $r \geq 2$ and $0<\theta<1$ are constants to control the subsample sizes (i.e., $h$  having larger   order than $h_2$) with  $C_1,C_2>0$.
A common pilot estimate $\widehat{V}_{\tilde{\ell}}$ is used for both $ \widehat{\mathrm{MSE}}_{h,n}(\ell)$ and    $\widehat{\mathrm{MSE}}_{h_2,n}(\ell)$.
We denote corresponding block estimates as $\widehat{\ell}_{opt,h}  $ and $\widehat{\ell}_{opt,h_2} $, and define an estimator of the target optimal block size $\ell_{opt,n}$ as
\[
\widehat{\ell}_{opt, n} \equiv     \left(  \frac{\widehat{\ell}_{opt, h}}{C_1^{\widehat{a}_n}}  \right)^r  \left(\frac{ h^{\widehat{a}_n}}{\widehat{\ell}_{opt, h}}\right)^{r-1}    \widehat{c}_n, \quad\quad\widehat{c}_n \equiv
r^{ \widehat{I}_n } \left(\frac{\log h_2 }{\log h} \right)^{ (r-1)\widehat{I}_n  (\log h)/\log(h/h_2) },
\]
where
\[
\widehat{a}_n \equiv   \frac{\log(\widehat{\ell}_{opt,h}/\widehat{\ell}_{opt,h_2})}{\log(h/h_2)},\quad
  \widehat{I}_n  \equiv \frac{1}{2}\left(\frac{\log (\widehat{\ell}_{opt, h}) - \widehat{a}_n \log h}{\log \log  h} + \frac{\log (\widehat{\ell}_{opt, h_2}) - \widehat{a}_n \log h_2}{\log \log h_2} \right)
\]
estimate the exponent $a\equiv a_{\alpha, m, m_2}\in(0,1)$ and indicator quantity $\tilde{I} \equiv -0.5  I(\alpha  = 0.5, m=1)+ I(\alpha m_2 = 1)$   
appearing in the Corollary~\ref{cor1} expression for $\ell_{opt,n}\approx K n^{a} (\log n )^{\tilde{I}}$.
 The estimator $\widehat{\ell}_{opt, n}$  has three components, where $ \widehat{\ell}_{h,opt}/C_1^{\widehat{a}_n}  $ estimates $ K h^{a} (\log h)^{\tilde{I}}$ while $\widehat{\ell}_{opt,h}/h^{\widehat{a}_n}$ captures
$ K (\log h)^{\tilde{I}}$ up to a  constant, and   $\widehat{c}_n$ is  scaling adjustment from $ \log n  \approx r \log h$.
The data-driven block estimator $\widehat{\ell}_{opt, n}$ is provably valid over differing forms for $ \ell_{opt, n}$ under LRD.

\begin{cor}
\label{cor4}  Let $h\equiv \lfloor C_1 n^{1/r}\rfloor$  and $h_2 =\lfloor C_2 n^{\theta/r} \rfloor$ (for $C_1,C_2>0$, $r\geq 2$, $\theta\in(0,1)$),
and suppose Corollary~\ref{cor3} assumptions hold.
 Then, as $n\to \infty$, the estimator $\widehat{\ell}_{opt, n}$ is consistent for $\ell_{opt,n}$   in that $\widehat{\ell}_{opt, n}/\ell_{opt,n} \stackrel{p}{\rightarrow} 1$ and, additionally,
\[
 \widehat{a}_n  \stackrel{p}{\rightarrow} a, \quad \widehat{I}_n \stackrel{p}{\rightarrow}\tilde{I}, \quad  \widehat{c}_n   \stackrel{p}{\rightarrow}  [r \theta^{(r-1)/(1-\theta)}]^{\tilde{I}}, \quad  \frac{\widehat{\ell}_{opt, h}}{ h^{\widehat{a}_n}} \cdot  \frac{1}{K (\log h)^{\tilde{I}}} \stackrel{p}{\rightarrow} \theta^{\tilde{I}/(1-\theta)},
\]
regarding constants $\tilde{I} \equiv -0.5  I(\alpha=0.5, m  = 1)+ I(\alpha m_2 = 1)$, $K\equiv K_{\alpha, m, m_2}>0$, and $a\equiv a_{\alpha, m, m_2}$
 for prescribing $\ell_{opt,n}\approx K n^{a} (\log n )^{\tilde{I}}$ under  Corollary~\ref{cor1}.
\end{cor}

We suggest a first subsample size $h = C_1  n^{1/2}$ of maximal possible order ($r=2$). We then take the pilot block to be $\tilde{\ell}=n^{1/2}$,   representing a reasonable choice under LRD and also satisfying  Theorem~\ref{cor3}-Corollary~\ref{cor4}  (i.e., $h/\tilde{\ell} + h\tilde{\ell}/n=O(1)$ then holds).
For a general rule in  numerical studies to follow,  we chose $C_1=9,C_2=12,\theta=0.95$ to keep subsamples adequately long under LRD.
We also consider a modified block estimation  rule
\begin{equation}
\label{eqn:rule}
\widehat{\ell}_n  = \min \{\lfloor n/20 \rfloor, \lfloor \widehat{\ell}_{opt, n}\rfloor \},
\end{equation}    to avoid overly large block selection in finite sample cases.  This variation retains consistency due to $\ell_{opt, n} = o(n)$ and  performs well over a variety of applications in Section~\ref{sec:6}.

\section{Extending the Scope of Statistics}
\label{sec:5}
Here we discuss extending  block selection and resampling variance estimation
to a larger class of statistics defined by functionals of empirical distributions.
Using a small notational change to develop this section, let us denote data from an observed time stretch as $Y_1=\tilde{G}(Z_1),\ldots, Y_n=\tilde{G}(Z_n)$  (rather than $X_t=G(Z_t)$) and let  $F_n \equiv n^{-1} \sum_{t=1}^n \delta_{Y_t}$ denote the corresponding
the empirical distribution, where $\delta_y$ denotes a unit point mass at $y\in\mathbb{R}$.
Consider a statistic
$T_n \equiv T(F_n)$, given by a real-valued functional $T(\cdot)$ of $F_n$,  which estimates a target parameter $T(F)$ defined by the process marginal distribution $F$.  A broad class of statistics and parameters can be expressed through such functionals, with some examples given below.\\

\noindent {\it Example 1:} Smooth functions $T_n$ of averages given by
\[T_n \equiv H\left( n^{-1}\sum_{t=1}^n \phi_1(Y_t),\ldots,  n^{-1}\sum_{t=1}^n \phi_l(Y_t) \right),
\]
involving a function $H:\mathbb{R}^l\rightarrow \mathbb{R}$  of $l\geq 1$ real-valued functions $\phi_j: \mathbb{R}\rightarrow \mathbb{R}$  for $j=1,\ldots,l$.
These statistics include ratios/differences of means as well as sample moments (\cite{L03}, ch.~5).  \\

\noindent {\it Example 2:} M-estimators $T_n$ defined as the solution to
\[ \frac{1}{n} \sum_{t=1}^n \psi(Y_t,T_n)=0
\]
for an estimating function with mean zero $\E [\psi(Y_t, T(F))]=0$.  This includes several types of location/scale or regression estimators
investigated in the LRD literature (cf.~\cite{B91, B10}).\\

\noindent {\it Example 3:} L-estimators $T_n$ defined through integrals as
\[ T_n = \int x J(F_n(x)) dF_n(x),
\]
involving a bounded function $J:[0,1]\rightarrow \mathbb{R}$.  These include
trimmed means $J(x) = {I}(\delta_1< x<\delta_2)/(\delta_2-\delta_1)$ (based on the indicator function and trimming proportions $\delta_1,\delta_2\in(0,1)$) along with Windsorized averages and Gini indices (cf.~\cite{S03}).  \\

For a fixed integer $k \geq 1$,   functionals defined by linear combinations  or products of components in  ``$k$-dimensional" marginal distributions might also be considered
(i.e., empirical distributions of $(Y_t,Y_{t+1},\ldots, Y_{t+k})$).
For simplicity, we use $k=1$.
Under regularity conditions \cite{F12, S03}, statistical functionals $T_n=T(F_n)$ as above are approximately linear and admit an expansion
\begin{equation}
\label{eqn:Tn}
T_n = T(F) + \frac{1}{n}\sum_{t=1}^n IF(Y_t,F)  + R_n,
\end{equation}
in terms of the influence function $IF(y, F)$,  defined as
\[
  IF(y, F) \equiv \lim_{\epsilon \downarrow 0}\frac{ T(   (1-\epsilon)F + \epsilon \delta_y)  - T(F) }{\epsilon},\quad y\in\mathbb{R},
\]
 and an  appropriately small remainder $R_n$; note that $\E [IF(Y_t,F)]=0$ holds. See \cite{DT89} and \cite{HM95}
  for such expansions with LRD Gaussian subordinated processes.

To link to our previous block resampling developments  (Sec.~\ref{sec:2}),
a statistic $T_n$ as in (\ref{eqn:Tn}) corresponds approximately  to an average $\bar{X}_n=\sum_{t=1}^n X_t/n$
of transformed LRD Gaussian observations
$X_1,\ldots,X_n$, where $X_t \equiv IF(Y_t,F) = IF(\tilde{G}(Z_t),F) = G(Z_t)$ has Hermite rank denoted by $m$ with $\alpha m<1$ here.   That is,
under appropriate conditions, the normalized statistic $n^{\alpha m/2} [T_n - T(F)] = n^{\alpha m/2}  \bar{X}_n   + o_p(1)  $ has a distributional limit
determined by   $\bar{X}_n$ (e.g.,~\cite{T75,T77}) with a limiting variance $\lim_{n\to \infty} n^{ \alpha m} \var(\bar{X}_n)  = v_{\infty, \alpha m}$ given by (\ref{xbar.var}) as before.  Results in \cite{BT19} also suggest that compositions
$X_t \equiv  G(Z_t) =IF(\tilde{G}(Z_t),F) $ may tend to produce Hermite ranks of $m=1$, in which case
$n^{\alpha /2} [T_n - T(F)]$ will be asymptotically normal with asymptotic variance  $v_{\infty, \alpha m }$.

To estimate   $v_{\infty, \alpha m}$ through  block resampling, we would ideally use  $X_1,\ldots,X_n$  to obtain  a variance estimator  as in Section~\ref{sec:2},  which we denote  as $\widehat{V}_{\ell,\alpha m} \equiv \widehat{V}_{\ell,\alpha m}(X)$. Then, all  estimation and  block properties from Sections~\ref{sec:3}-\ref{sec:4} would apply.  Unfortunately,   $F$ is generally unknown in practice so that $\{X_t\equiv IF(Y_t,F)\}_{t=1}^n$ are unobservable from the  data $Y_1,\ldots,Y_n$.  Consequently, $\widehat{V}_{\ell,\alpha m }(X)$ represents an oracle  estimator.
In Sections~\ref{sec5.1}-\ref{sec5.2}, we detail two block-based strategies for estimating    $v_{\infty, \alpha m}$
based on either a substitution   method or block jackknife.  In both cases, these approaches
can be as good as the oracle estimator $\widehat{V}_{\ell,\alpha m}(X)$ under some conditions.
These resampling results under LRD have counterparts to the SRD case \cite{K89, PP02}, though we non-trivially include L-estimation in addition to M-estimation.

\subsection{Substitution method}
\label{sec5.1}

Classical substitution (i.e., plug-in)  estimates $F$ in the influence function $IF(y,F)$ with its empirical version $F_n$ (cf.~\cite{ PP02, S03}) and   develops observations as $\widehat{X}_1, \ldots, \widehat{X}_n$ with $\widehat{X}_t \equiv IF(Y_t,F_n)$.   For example, in a smooth function $T_n = H(\bar{Y}_n)$ of the average $\bar{Y}_n$, we have $IF(y, F_n) =H^\prime(\bar{Y}_n)(y- \bar{Y}_n)$, where $H^{\prime}$ denotes the derivative of $H$.
We denote a resampling variance estimator computed from such observations as  $ \widehat{V}_{\ell,\alpha m}( \widehat{X})$.

To  compare $ \widehat{V}_{\ell,\alpha m}( \widehat{X})$  to the oracle estimator $\widehat{V}_{\ell,\alpha m}(X)$,
we require bounds between estimated $IF(y,F_n)$ and true influence functions $IF(y,F)$.  For weakly dependent processes and M-estimators, \cite{K89} considered pointwise expansions of $IF(y,F)-IF(y,F_n)$ as linear combinations of other functions in $y$.  We need to generalize the concept of such expansions to accommodate LRD and more general functionals (e.g., L-estimators) as follows.\\

\noindent \textbf{Condition-I:}  There exist random variables $U_{1,n},U_{2,n},W_n$ and real constants $c,d\in\mathbb{R}$, $C>0$  such that, for any generic real values  $y_1,\ldots,y_k$ and $k \geq 1$, it holds that
\[
  \Bigg| \frac{1}{k} \sum_{j=1}^k IF(y_j, F )  -    \frac{1}{k}\sum_{j=1}^k IF(y_j, F_n )  +  U_{1,n}\Bigg| \leq |W_n| + |U_{2,n}|\left(\int_{c}^{d}  \Bigg| \frac{1}{k}\sum_{j=1}^k h_\lambda (y_j)   \Bigg|^2 d \lambda\right)^{1/2},
\]
 where $|W_n|= O_p( n^{-\alpha m})$;   $ |U_{1,n}|, |U_{2,n}|=O_p(n^{-\alpha m/2})$;  and, as indexed by $\lambda \in[c,d]$, $h_\lambda(\cdot)$ denotes a real-valued function  such that $h_\lambda (Y_t)=h_\lambda(\tilde{G}(Z_t))$ has mean zero, variance $\E[h_\lambda (Y_t)]^2 \leq C$, and Hermite rank of at least $m$ (the  rank of $G(Z_t) = IF(Y_t,F)$).\\

For context, if we set $\alpha m=1$ above and skip the notion of Hermite rank, then Condition-I would include, as a special case, an assumption  used by  \cite{K89} with weakly dependent processes. However, under LRD, we need to explicitly incorporate Hermite ranks in bounds.  If we define  $m_y \geq 1$ as
the Hermite rank of an indicator function $I(Y_t \leq y) = I(\tilde{G}(Z_t)\leq y)$ for $y\in \mathbb{R}$, then the smallest rank $m^* \equiv  \{m_y:y\in\mathbb{R}\}$ is known to be useful for describing convergence of the empirical distribution $[F_n(\cdot) - F(\cdot)]$ (cf.~\cite{DT89}).  One general way to ensure any function $h_\lambda(\cdot)$ appearing in Condition-I has Hermite rank of at least $m$ (the rank of $IF(Y_t,F)$)
is that $m=m^*$.  The reason is that $m^*$ sets a lower bound on the Hermite rank of any function of $Y_t$ (cf.~(2.5) of \cite{DT89}).  Such equality $m=m^*$ appears implicit in work of \cite{HM95} on statistical functionals under LRD
and holds automatically when $m=1$.  We show next that  the statistics $T_n$ in Examples~1-3 can satisfy Condition-I.

\begin{theorem}
\label{theorem4} For $Y_t = \tilde{G}(Z_t)$,  suppose $X_t \equiv G(Z_t) = IF(Y_t,F)$ has Hermite rank $m\geq 1$ with $\alpha m <1$. Then, Condition-I holds if the functional $T_n$ is as in \\[.1cm]
\textrm{(i)}  Example 1 (smooth function)  where  $\phi_1,\ldots, \phi_l$ are bounded functions; first partial derivatives of $H :\mathbb{R}^l\rightarrow \mathrm{R}$  are Lipschitz in a neighborhood of $(\E[\phi_1(Y_t)], \ldots, \E[\phi_l(Y_t)])$; and  either $m=m^*$ holds or $m=\min\{\mbox{Hermite rank of $\phi_j(Y_t): 1 \leq j \leq l$}\}$.
 \\[.1cm]
\textrm{(ii)}  Example 2 (M-estimation) where  a constant $C>0$ and a neighborhood $N_0$ of $T(F)$ exist such that
$|\psi(y,\theta)|\leq C$ on $\mathbb{R}\times N_0$;   $\dot{\psi}\equiv \partial \psi/\partial \theta$  exists and $|\dot{\psi}(y,\theta)|\leq C$ on $\mathbb{R}\times N_0$;  $|\dot{\psi}(y,\theta_1) - \dot{\psi}(y,\theta_2)|\leq C|\theta_1 - \theta_2|$ for $y\in\mathbb{R}$, $\theta_1,\theta_2\in N_0$; $\E \dot{\psi}(Y_t,T(F)) \neq 0 $; and either $m=m^*$ holds or the Hermite rank of $\psi(Y_t, \theta)$ remains the same for $\theta\in N_0$.  \\[.1cm]
\textrm{(iii)}  Example 3 (L-estimation) where $J$ is bounded and Lipschtiz on $[0,1]$ with $J(t)=0$ when $t\in [0,\delta_1] \cup [\delta_2,1]$  for some $0 < \delta_1<\delta_2 < 1$; and either $m=m^*$ holds or     $m \leq \min\{m_y:  y_1\leq y\leq  y_2\}$ for some
real $y_1<y_2$ with $0<F(y_1)<\delta_1 < \delta_2<F(y_2)<1$.
\end{theorem}

Theorem~\ref{theorem4} assumptions for Examples 1-2, dropping Hermite rank conditions,  match those of \cite{K89}.
 Smooth function statistics in Example~1 have influence functions $X_t=IF(Y_t,F)$ as a linear combination of the baseline functions $\phi_j(Y_t)$, $1\leq j \leq l$,
 so that the smallest  Hermite rank among these typically gives the  Hermite rank $m$ of $IF(Y_t,F)$.
 In   M-estimation, the Hermite rank   of $X_t \equiv G(Z_t) = IF(Y_t,F)$ matches that of $\psi(Y_t, T(F))$ and it is sufficient that $\psi(Y_t, \theta)$ maintains the same rank $m$ in a $\theta$-neighborhood of $T(F)$; the latter condition is mild and  implies that the rank of  $ \dot{\psi}(Y_t, T(F))$ must be at least $m$, which is important as $\dot{\psi}(\cdot, T(F))$ arises in Condition-I under M-estimation.  To illustrate  with a standard normal $Y_t=Z_t$, M-estimation of the process mean uses $\psi(Z_t,\theta)=Z_t-\theta$ with a constant Hermite rank of 1 as a function of $\theta$ and a derivative $\dot{\psi}(Y_t,T(F))=-1$
of infinite rank; similarly, Huber-estimation uses  $\psi(Z_t,\theta)= \max\{-c ,\min\{Z_t-\theta,c\} \}$ (for some $c>0$) which has constant rank 1 for $\theta$ in a neighborhood of $T(F)=0$ here, while the derivative
$\dot{\psi}(Z_t,T(F))=  I(|Z_t|\leq c) $ has rank 2.
For general L-estimation, conditions on the Hermite ranks $m_y$ of
indicator functions $I(Y_t \leq y)$ (or the empirical distribution $F_n(y)$) are necessary, particularly when trimming percentages $\delta_1,\delta_2$ are involved; in this case, we may use the rank $m_y$ of $F_n(y)$ over
a $y$-region ($[y_1,y_2]$) that is not trimmed away.

Theorem~\ref{theorem5}  establishes that  the oracle  resampling estimator $\widehat{V}_{\ell,\alpha m}(X)$ (true influence) and the plug-in version  $\widehat{V}_{\ell,\alpha m}(\widehat{X})$ (estimated influence) are often close to the extent that  the latter is as good as the former. Blocks can be either OL/NOL below.

\begin{theorem}
   \label{theorem5}
For $Y_t = \tilde{G}(Z_t)$,  suppose $X_t \equiv G(Z_t) = IF(Y_t,F)$ has Hermite rank $m\geq 1$ with $\alpha m <1$, Condition-I holds, and   $\ell^{-1} + \ell/n\rightarrow 0$ as $n\to \infty$.  Then,
\[
  \widehat{V}_{\ell,\alpha m}(\widehat{X}) = \widehat{V}_{\ell,\alpha m}(X)  +  O_p( (\ell/n)^{\alpha m /2} ) + O_p(n^{-\alpha m /2}).
\]
\end{theorem}

Theorem~\ref{theorem5} is the LRD analog of a result by \cite{K89} for weakly dependent processes (i.e., setting $\alpha m =1$ above).  As in the SRD case, the difference between estimators
is often  no larger than the estimation error $O_p( (\ell/n)^{\min\{\alpha,1/2\}})$ from the standard deviation of the oracle  $\widehat{V}_{\ell,\alpha m}(X)$  (Theorem~\ref{theorem2}).   Consequently,  optimal block   orders and convergence rates for $\hat{V}_{\ell,\alpha m}(X)$ (Section~\ref{sec:4.1}) generally apply to the substitution version $\widehat{V}_{\ell,\alpha m}(\widehat{X})$.  The block rule of Section~\ref{sec:4.2} can also be applied to $\widehat{X}_1,\ldots,\widehat{X}_n$, which we illustrate in Section~\ref{sec:6}.

\subsection{Block jackknife (BJK) method}
\label{sec5.2}
For estimating the asymptotic variance $ v_{\infty, \alpha m}$ of the functional $T_n$, a block jackknife (BJK) estimator is possible under LRD.  BJK uses only OL data blocks, as NOL blocks are generally invalid (Remark 4.1, \cite{K89}) .     For  $j=1,\ldots,N\equiv n-\ell +1$, we compute
the functional $T_n^{(j)}$ after removing observations in  $j$th OL block $(Y_j,\ldots,Y_{j+\ell-1})$  from the data $(Y_1,\ldots,Y_n)$.   The BJK estimator of $ v_{\infty, \alpha m}$ is then
\[
\widehat{V}_{\ell,\alpha m,\OL}^{{\scriptscriptstyle BJK}}   = \frac{(N-1)^2}{\ell^2}  \frac{\ell^{\alpha m}}{N}\sum_{j=1}^N (T_n^{(j)} - \bar{T}_n)^2,\quad \bar{T}_n \equiv  \frac{1}{N} \sum_{j=1}^N T_n^{(j)}.
\]
   Unlike the plug-in method (Sec.~\ref{sec5.1}), BJK does not  involve influence functions, but uses repeated evaluations  of the functional. For the sample mean statistic $T_n= \sum_{t=1}^n Y_t/n$, the BJK estimator     matches the
plug-in estimator $\widehat{V}_{\ell,\alpha m}(\widehat{X}) \equiv \widehat{V}_{\ell,\alpha m, \OL}(\widehat{X})$   with OL blocks (cf.~\cite{K89}).
More generally, these two estimators may differ, though not substantially, as shown in Theorem~\ref{theorem6}.
To state the result, for each OL data block $j=1,\ldots,N$, we define a remainder   $S_n^{(j)} \equiv T_n^{(j)} - T_n - M_n^{(j)}$, due to a type of Taylor expansion of $T_n^{(j)}$ about $T_n$, where
\begin{eqnarray}
\nonumber M_n^{(j)} \equiv \frac{1}{n-\ell} \sum_{\substack{1 \le  t \le n, \\ t \not\in [j, j+\ell-1]}} \widehat{X}_{t} - \frac{1}{n} \sum_{t=1}^n \widehat{X}_{t}
\end{eqnarray}
involves an average of  estimated values $\{\widehat{X}_i \equiv IF(Y_i,F_n)\}_{i=1}^n$  after removing  the $j$th block.

\begin{theorem}
   \label{theorem6}
For $Y_t = \tilde{G}(Z_t)$,  suppose $X_t \equiv G(Z_t) = IF(Y_t,F)$ has Hermite rank $m\geq 1$ with $\alpha m <1$, and
that the OL block plug-in estimator $\widehat{V}_{\ell,\alpha m, \OL}(\widehat{X})$ is consistent. Then,
\[
   \widehat{V}_{\ell,\alpha m}^{{\scriptscriptstyle BJK}}  = \widehat{V}_{\ell,\alpha m, \OL}(\widehat{X})+ O_p(\ell/n)
\]
holds as $n\to\infty$ if $ \ell^{\alpha m} \sum_{j=1}^N [S_n^{(j)}]^2/N = O_p( \ell^4/[n^2 (N-1)^2] )$; the latter is true under Theorem~\ref{theorem4} assumptions for Examples~1-3.
\end{theorem}
The above difference $O_p(\ell/n)$  between BJK and plug-in estimators holds similarly under weak dependence (akin to setting $\alpha m=1$ above), which improves the bound $O_p(\ell^{3/2}/n)$ originally given by \cite{K89} (Theorem~4.2).  Theorems~\ref{theorem5}-\ref{theorem6}  show  that BJK can also differ
no more from the oracle estimator  $\widehat{V}_{\ell,\alpha m, \OL}(X)$ than the plug-in estimator $\widehat{V}_{\ell,\alpha m, \OL}(\widehat{X})$.

\section{Numerical Illustrations and Applications}
 
\label{sec:6}
\subsection{Illustration of MSE over block sizes}
\label{sec6.1}
\begin{figure}[ht]
\includegraphics[scale=0.5]{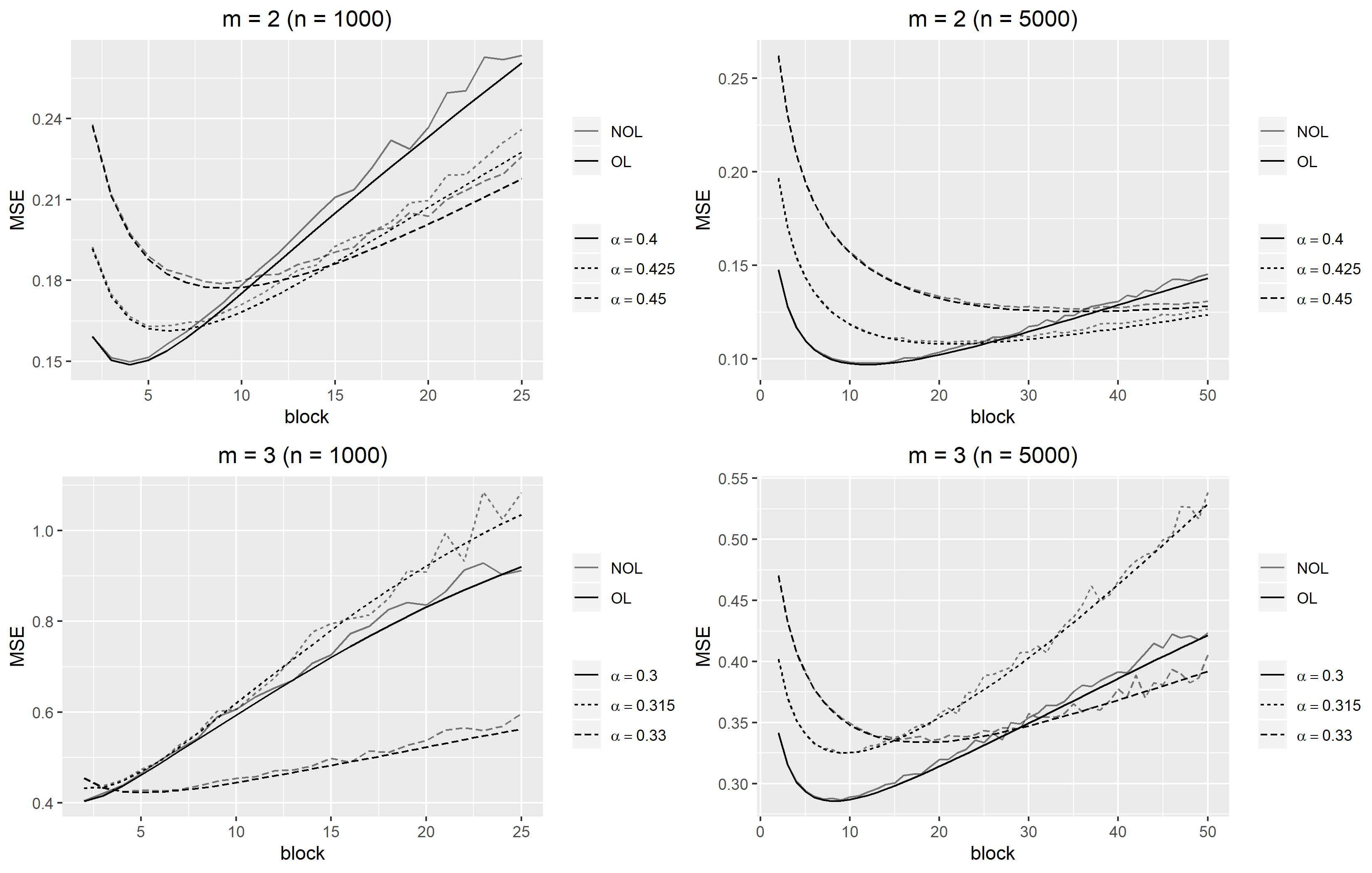}
\caption{MSE curves (over block length $\ell$) for resampling estimators  with different LRD processes $(m,\alpha)$.}
\label{fig:Optimal block size}
\end{figure}

Here we describe an initial numerical study of the MSE-behavior of   resampling variance estimators under LRD.
In particular, results of Section~\ref{sec:3} suggest that OL/NOL resampling blocks should induce identical large-sample performances under strong dependence (e.g., $\alpha<1/2$)
 and that optimal blocks should   generally {\it decrease}   in size   as the covariance strength
  increases   (cf.~Sec~\ref{sec:4.1}).
LRD series were   generated    as
$X_t = H_2(Z_t)$ or $X_t = H_3(Z_t)$, using three values of the memory parameter with $\alpha < 1/m$ for $m=2$ or $m=3$,
based on a  standardized Fractional Gaussian process $Z_t$   with covariances
 as in (\ref{eqn:zcov})  (i.e., $H=(2-\alpha)/2$).
For each simulated series, OL/NOL block-based estimators $\widehat{V}_{\ell,\am}$ of the variance $v_{n,\am}=n^{\am} \var(\bar{X}_n)$ were computed over a sequence of block sizes $\ell$.  Repeating this procedure over 3000 simulation runs and
averaging  differences $(\widehat{V}_{\ell,\am}-v_{n,\am})^2/ v_{n,\am}^2$ produced approximations
of   standardized  MSE-curves $\E(\widehat{V}_{\ell,\am}-v_{n,\am})^2/ v_{n,\am}^2$, as   shown in Figure~\ref{fig:Optimal block size}
with sample sizes $n=1000$ or $5000$.  The MSE curves are quite close between OL/NOL blocks, particularly as sample sizes increase to $n=5000$,
in agreement with theory.
  Also, as suggested by Section~\ref{sec:4.1}, MSEs should improve at the best block choice as    covariance strength
  increases under LRD ($\alpha \downarrow$), which is visible in  Figure~\ref{fig:Optimal block size}.    Table~\ref{table:Optimal block size} presents    best block lengths from the figure, showing that  optimal blocks decrease
 for these  LRD processes with decreasing $\alpha$. The supplement \cite{ZLN22} provides additional simulation studies to further illustrate   bias/variance
 behavior of resampling estimators.

\begin{table}[ht]
\centering
\caption{Optimal OL block sizes for LRD series with Hermite ranks $m = 2$ and $3$}
\begin{tabular}{lccccccc}
  & \multicolumn{3}{c}{$m=2$} & &\multicolumn{3}{c}{$m=3$} \\
\multicolumn{1}{r}{$\alpha$} &  0.400  &  0.425 &  0.450   & &  0.300  & 0.315  &  0.330 \\ \cline{2-4} \cline{6-8}
$n=1000$  & 4 &6 &9 &&  2 &3& 5\\
$n=5000$  & 12 &21 &35&& 8 &9 &18\\
\hline
\end{tabular}
\label{table:Optimal block size}
\end{table}

\subsection{Resampling variance estimation by empirical  block size}
\label{subsec:msevar}
We next examine empirical block choices for  resampling  variance estimation of the sample mean and provide comparison to other  approaches under LRD.   Application to another functional   is then considered.

We use  the  data-based rule (\ref{eqn:rule}) of Section~\ref{sec:4.3} for choosing a block size.
  We first compare  resampling estimators $\widehat{V}_{\ell}$  of the sample mean's variance $v_{n,\am}$ between block selections $\ell=\widehat{\ell}_{n}$ and $\ell=\lfloor n^{1/2} \rfloor$, where the latter  represents a reasonable choice under LRD by theory in Section~\ref{sec:4.1}.
    OL blocks are used along with local Whittle estimation $\widehat{\alpha m}$  of the memory parameter $\am$ (Sec.~\ref{sec:4.2}).
    Similarly to Section~\ref{sec6.1}, we simulated samples from LRD processes defined by $X_t = H_m(Z_t)$ for $m = 2$ with $\alpha = 0.20, 0.45$ or $m = 3, \alpha = 0.20, 0.30$
and approximated the
 MSE $\E(\widehat{V}_{\ell} -v_{n,\am})^2/v_{n,\am}^2$   using $500$ simulations.
 Table~\ref{tab:blockmse} provides these results.  Estimated blocks $\widehat{\ell}_n$ are generally better than the default   $\lfloor n^{1/2} \rfloor$,
 though the latter is also competitive. The default seems preferable with a small sample size and particularly strong dependence (e.g., $n=500$, $m = 3, \alpha = 0.2$),
 but
 empirical block selections
  show
  improved MSEs
   with  increased sample sizes $n=1000, 2000$ under LRD.

For comparison against   resampling variance estimators, we also consider the Bartlett-kernel heteroskedasticity  and autocorrelation consistent (HAC) estimator \cite{W80} and the memory and autocorrelation consistent (MAC) estimator \cite{R05}, whose large-sample properties have been  studied for the sample mean with linear LRD
processes (cf.~\cite{ADG09, GRS99}), but not  for transformed LRD series $X_t=G(Z_t)$.
As numerical suggestions from \cite{ADG09},    we implemented HAC and MAC estimators of the sample mean's variance using
 bandwidths  $\lfloor n^{1/5} \rfloor$,   $\lfloor n^{4/5} \rfloor$, respectively; the HAC approach further used  local Whittle estimation of the memory parameter $\am$, like the resampling estimator.  The MSEs of HAC/MAC estimators are given in
 Table~\ref{tab:hacmacmse} (approximated from 500 simulation runs) for comparison against the resampling estimators in Table~\ref{tab:blockmse} with the same processes.
For the process $X_t = H_2(Z_t)$ with $\alpha = 0.45$,  HAC/MAC estimators emerge as slightly better than the resampling approach with estimated block sizes,
though the  resampling estimator  outperforms  HAC/MAC estimators as the dependence increases (smaller $\alpha$) or as the Hermite rank increases ($m=3$).  With small sample sizes $n$ and strong dependence, the HAC estimator  can exhibit large MSEs, indicating that the  bandwidth $\lfloor n^{1/5} \rfloor$ is  perhaps too small
for the non-linear LRD series in these settings.
In comparison, the empirical block  selections with resampling  estimators show    consistently reasonable  MSE-performance among all cases, which is appealing.

We further consider a different statistical functional with resampling   estimators and  empirical blocks  $\widehat{\ell}_{n}$.
In the notation of Section~\ref{sec:5}, we simulated stretches $Y_1,\ldots,Y_n$ of LRD processes defined by $Y_t = H_2(Z_t)$ or $Y_t= \sin(Z_t)$
and considered an L-estimator $T_n$ as a 40\% trimmed mean  based on the empirical distribution $F_n$  (i.e., $\delta_1 =1-\delta_2=0.2$ in Example~3, Sec.~\ref{sec:5}).  For either process, the influence function $X_t\equiv IF(Y_t,F) = Y_t I( F^{-1}(0.2) < Y_t < F^{-1}(0.8))/0.6$ has Hermite rank $m=1$, where $F$ and $F^{-1}$ denote the distribution and quantile functions, respectively, of $Y_t$.
 To estimate
the variance, say $v_{n,\am}$, of $n^{\alpha m/2} T_n$, we apply the substitution method (Sec.~\ref{sec5.1}).  That is, using estimated influences
$\widehat{X}_t \equiv IF(Y_n,F_n)= Y_t I( F_n^{-1}(0.2) < Y_t < F_n^{-1}(0.8))/0.6$, we obtain  an estimator  $\widehat{\alpha m}$ of the memory-parameter by  local Whittle estimation and  compute  a plug-in resampling variance estimator $\widehat{V}_{\ell}(\widehat{X})$.    Table~\ref{tab:trimmedblockmse} provides   MSEs  (i.e., $\E(\widehat{V}_{\ell}(\widehat{X}) - v_{n,\am})^2/v_{n,\am}^2$    approximated from 500 simulation runs)
with block choices $\ell=\widehat{\ell}_n$ or $\lfloor n^{1/2}\rfloor$ over sample sizes $n=500, 1000, 2000$.  The empirical block selections perform  better than the default  $\lfloor n^{1/2}\rfloor$ with
the plug-in   variance estimator here, though the choice $\lfloor n^{1/2}\rfloor$ appears also reasonable.

\begin{table}[t]
\centering
\caption{MSE for resampling  estimators     $\widehat{V}_{\ell}$  of  sample mean's variance based on blocks $\ell=\widehat{\ell}_n$ or $\ell= \lfloor n^{1/2}\rfloor$ with  LRD series $X_t \equiv G(Z_t)$ with $G= H_2(Z_t)$ $(m=2)$ and $G(Z_t) = H_3(Z_t)$ $(m=3)$.}

\label{tab:blockmse}
\begin{tabular}{cccccccccccc}
\\
\cline{1-12}
& &  &    & \multicolumn{2}{c}{$n=500$} & & \multicolumn{2}{c}{$n=1000$} & & \multicolumn{2}{c}{$n=2000$}\\\cline{5-12} \vspace*{.1cm}
\raisebox{-.1cm}{$G/m$}    & &  \raisebox{-.1cm}{$\alpha$}   &  &    \raisebox{-.1cm}{$ \widehat{\ell}_{n}$}   &  \raisebox{-.1cm}{$ \lfloor n^{1/2}\rfloor$}   & &   \raisebox{-.1cm}{$ \widehat{\ell}_{n}$}  &  \raisebox{-.1cm}{$ \lfloor n^{1/2}\rfloor$}   & &   \raisebox{-.1cm}{$ \widehat{\ell}_{n}$} & \raisebox{-.1cm}{$ \lfloor n^{1/2}\rfloor$ } \\
\cline{1-12}
\multirow{2}{*}{$m=2$} & & $0.20$  & & 0.294 & 0.316  & & 0.236 &  0.248 & & 0.180 & 0.214                                 \\
& & $0.45$  & & 0.270 & 0.312  & & 0.269 &  0.293 & & 0.280 & 0.292         \\
\cline{1-12}
\multirow{2}{*}{$m=3$} & & $0.20$   & & 0.917 & 0.754  & & 0.805 &  0.858 & & 0.455 & 0.495       \\
& & $0.30$  & & 0.270 & 0.294  & & 0.396 &  0.413 & & 0.393 & 0.393       \\
\cline{1-12}
\end{tabular}
\end{table}

\begin{table}[t]
\centering
\caption{MSE for  HAC  and MAC  estimators of sample mean's variance with  LRD series $X_t \equiv G(Z_t)$ with $G= H_2(Z_t)$ $(m=2)$ and $G(Z_t) = H_3(Z_t)$ $(m=3)$.}

\label{tab:hacmacmse}
\begin{tabular}{llllllllllll}  \\
\cline{1-12}
& &  &    & \multicolumn{2}{c}{$n=500$} & & \multicolumn{2}{c}{$n=1000$} & & \multicolumn{2}{c}{$n=2000$}   \\
\cline{5-12}
$G/m$   & & $\alpha$  & & HAC & MAC   & & HAC & MAC   & &  HAC & MAC   \\ \cline{1-12}
\multirow{2}{*}{$ m=2 $} & & $0.20$  & & 8.297 & 0.336  & & 0.989 &  0.250 & & 0.286 & 0.174     \\
& & $0.45$  & & 0.256 & 0.252  & & 0.265 &  0.257 & & 0.279 & 0.267   \\
\cline{1-12}
\multirow{2}{*}{$ m=3 $} & & $0.20$   & & 26.43 & 1.667  & & 1.492 &  0.797 & & 1.020 & 0.804     \\
& & $0.30$  & &  0.515 & 0.442  & & 0.379 &  0.378 & & 0.369 & 0.358    \\
\cline{1-12}
\end{tabular}
\end{table}

\begin{table}[t]
\centering
\caption{MSE for plug-in resampling  estimators     $\widehat{V}_{\ell}(\widehat{X})$  of the variance of the 40\% trimmed mean  statistic $T_n$
based on blocks $\ell=\widehat{\ell}_n$ or $\ell= \lfloor n^{1/2}\rfloor$ with  LRD series $Y_t \equiv  H_2(Z_t)$ or $Y_t = \sin(Z_t)$.}

\label{tab:trimmedblockmse}
\begin{tabular}{llllllllllll}
\\
\cline{1-12}
& &  &    & \multicolumn{2}{c}{$n=500$} & & \multicolumn{2}{c}{$n=1000$} & & \multicolumn{2}{c}{$n=2000$}                                                       \\ \cline{5-12}
\vspace*{.1cm}
\raisebox{-.1cm}{$Y_t$}    & &  \raisebox{-.1cm}{$\alpha$}   &  &    \raisebox{-.1cm}{$ \widehat{\ell}_{n}$}   &  \raisebox{-.1cm}{$ \lfloor n^{1/2}\rfloor$}   & &   \raisebox{-.1cm}{$ \widehat{\ell}_{n}$}  &  \raisebox{-.1cm}{$ \lfloor n^{1/2}\rfloor$}   & &   \raisebox{-.1cm}{$ \widehat{\ell}_{n}$} & \raisebox{-.1cm}{$ \lfloor n^{1/2}\rfloor$ } \\ \hline
\multirow{2}{*}{$H_2(Z_t)$} & & $0.20$  & & 0.463 & 0.494  & & 0.383 &  0.413 & & 0.366 & 0.402                                \\
& & $0.45$  & & 0.568 & 0.603  & & 0.565 &  0.585 & & 0.565 & 0.575         \\\cline{1-12}
\multirow{2}{*}{$\sin(Z_t)$} & & $0.20$   & & 0.608 & 0.641  & & 0.597 &  0.613 & & 0.590 & 0.606       \\
& & $0.30$  & & 0.561 & 0.589  & & 0.565 &  0.582 & & 0.550 & 0.565      \\\cline{1-12}
\end{tabular}
\end{table}

\subsection{Resampling distribution estimation by  empirical  block size}
Block selection also plays an important role in other resampling inference,   such as approximating full sampling distributions with block bootstrap for purposes of tests and
confidence intervals.  While optimal block sizes for distribution estimation are difficult and unknown under LRD,
we may apply blocking notions developed here for  guidance.  For distributional approximations of sample means  and other statistics as in Section~\ref{sec:5},
the block bootstrap is valid with transformed LRD series when a normal limit exists  (e.g., Hermite rank  $m=1$)  \cite{L93}.  Such normality may occur commonly  in practice  \cite{BT19}  and can be further assessed  as described in Section~\ref{sec6.4}.

To study empirical blocks  for distribution estimation with the bootstrap,
 we consider two LRD processes as $Y_t =  \sin(Z_t)$ or $Y_t =  Z_t + 20^{-1}H_2(Z_t)$ defined by Gaussian $\{Z_t\}$ as before with
 memory exponent  $\alpha$.  Based on a size $n$ sample,  block bootstrap is applied to approximate the distribution of
  $ \Delta_n\equiv n^{\widehat{\alpha m}/2} |T_n- \theta|$, where $T_n\equiv T(F_n)$ represents either the sample mean or the 40\% trimmed mean (Example~3, Sec.~\ref{sec:5})),
  while $\theta\equiv T(F)$ denotes the corresponding process mean or trimmed mean parameter.
    In sample mean case,
  we compute $\widehat{\alpha m}$ using local Whittle estimation with data stretch $X_1=Y_1,\ldots,X_n=Y_n$ and define a bootstrap
  average $\bar{X}_{b\ell}^{*}$ by resampling $b\equiv \lfloor n/\ell \rfloor$ OL data blocks of length $\ell$ (see~Sec~\ref{sec2.2});
   the bootstrap version of $\Delta_n$ is then  $ \Delta_n^* \equiv b^{1/2} \ell^{\widehat{\alpha m}/2} |\bar{X}^*_{b\ell} - \E_* \bar{X}^*_{b\ell}|$ (cf.~\cite{L93}), where
   $\E_* \bar{X}^*_{b\ell} = \sum_{i=1}^{n-\ell+1} \bar{X}_{i,\ell}/(n-\ell+1)$ is a bootstrap expected average.  In the trimmed mean case,
   the estimator $\widehat{\alpha m}$ and the bootstrap approximation $ \Delta_n^*$ are similarly defined from estimated values
   $\widehat{X}_t \equiv IF(Y_n,F_n)= Y_t I( F_n^{-1}(0.2) < Y_t < F_n^{-1}(0.8))/0.6$, $t=1,\ldots,n$.
   We construct 95\% bootstrap confidence intervals for $\theta$   by approximating the 95th percentile of $\Delta_n$
   with the bootstrap counterpart from  $\Delta_n^*$  (based on 200 bootstrap re-creations).  Note  that, for these processes and statistics, the effective Hermite rank is $m=1$ (i.e., the rank of $X_t = IF(Y_t,F)$) so that the bootstrap should be valid in theory.

 We used the empirical rule (\ref{eqn:rule}) as a guide for selecting a block length $\ell$.  Table~\ref{tab:MBBCI}  shows the empirical coverages of 95\% bootstrap intervals with samples of size $n=1000$ or $n=5000$ (based on 500 simulation runs).  For strongest LRD   $\alpha=0.02$,  bootstrap intervals exhibit under-coverage, as perhaps expected, though accuracy improves
 with increasing sample size $n$ in this case.  The bootstrap performs well in the other cases of long-memory.
 The coverage rates of bootstrap intervals are closer to the nominal
 level with empirically chosen blocks $\widehat{\ell}_n$ compared to a standard choice
 $\lfloor n^{1/2} \rfloor$, for both the sample mean and trimmed mean.  This    suggests that the driven-data rule for blocks provides a reasonable guidepost for resampling distribution estimation, as an application beyond  variance estimation.

\begin{table}[t]
\centering
\caption{Empirical coverage probabilities of 95\% block bootstrap confidence intervals for the process mean or 40\% trimmed mean, based on block sizes
  $\ell = \widehat{\ell}_{n}$ or $\ell=\lfloor n^{0.5}\rfloor$, with LRD series as $ \sin(Z_t)$  or  $Z_t + 20^{-1} H_2(Z_t)$.}
\label{tab:MBBCI}
\begin{tabular}{cccccccccccc}
\multicolumn{12}{c}{Mean}\\
\cline{1-12}
& &  &    & \multicolumn{2}{c}{$\alpha =0.20$} & & \multicolumn{2}{c}{$\alpha =0.50$} & & \multicolumn{2}{c}{$\alpha =0.80$}                                                       \\ \cline{5-12}
              & & $\ell$  & & $n=1000$ & $ 5000$   & & $n=1000 $& $ 5000$   & & $n=1000 $& $ 5000$                           \\ \cline{3-12}\vspace*{.1cm}
\multirow{2}{*}{\raisebox{-.1cm}{$ \sin(Z_t)$}} & & \raisebox{-.1cm}{$\widehat{\ell}_{n}$}  & & \raisebox{-.1cm}{0.820} & \raisebox{-.1cm}{0.866}  & & \raisebox{-.1cm}{0.932} &  \raisebox{-.1cm}{0.944} & & \raisebox{-.1cm}{0.964} & \raisebox{-.1cm}{0.958}                                 \\
& & $\lfloor n^{1/2} \rfloor$  & & 0.766 & 0.814  & & 0.878 &  0.914 & & 0.958 & 0.934         \\\cline{1-12}\vspace*{.1cm}
\multirow{2}{*}{\raisebox{-.1cm}{$Z_t + \frac{1}{20} H_2(Z_t)$}} & & \raisebox{-.1cm}{$\widehat{\ell}_{n}$}   & &\raisebox{-.1cm}{0.848} & \raisebox{-.1cm}{0.864}  & & \raisebox{-.1cm}{0.922} &  \raisebox{-.1cm}{0.930} & & \raisebox{-.1cm}{0.942} & \raisebox{-.1cm}{0.960} \\
& & $\lfloor n^{1/2} \rfloor$ & & 0.806 & 0.854  & & 0.914 &  0.894 & & 0.938 & 0.952          \\ \cline{3-12}
 & &   & & &    & & &   & &  &  \\
\multicolumn{12}{c}{Trimmed Mean}\\
\cline{1-12}
& &  &    & \multicolumn{2}{c}{$\alpha =0.20$} & & \multicolumn{2}{c}{$\alpha =0.50$} & & \multicolumn{2}{c}{$\alpha =0.80$}                                                       \\ \cline{5-12}
                 & & $\ell$  & & $n=1000$ & $ 5000$   & & $n=1000 $& $ 5000$   & & $n=1000 $& $ 5000$                           \\ \cline{3-12}\vspace*{.1cm}
\multirow{2}{*}{\raisebox{-.1cm}{$ \sin(Z_t)$}} & & \raisebox{-.1cm}{$\widehat{\ell}_{n}$ } & &\raisebox{-.1cm}{0.798} &\raisebox{-.1cm}{0.848}  & &\raisebox{-.1cm}{0.934} &  \raisebox{-.1cm}{0.942} & &\raisebox{-.1cm}{0.964}& \raisebox{-.1cm}{0.950}                                 \\
& & $\lfloor n^{1/2} \rfloor$  & & 0.752 & 0.792 & & 0.892 &  0.920 & & 0.962 & 0.946         \\\cline{1-12}\vspace*{.1cm}
\multirow{2}{*}{\raisebox{-.1cm}{$Z_t + \frac{1}{20} H_2(Z_t)$}} & & \raisebox{-.1cm}{$\widehat{\ell}_{n}$}   & & \raisebox{-.1cm}{0.864} & \raisebox{-.1cm}{0.876}  & & \raisebox{-.1cm}{0.930} & \raisebox{-.1cm}{0.916}& & \raisebox{-.1cm}{0.956} & \raisebox{-.1cm}{0.980}  \\
& & $\lfloor n^{1/2} \rfloor$ & & 0.828 & 0.852  & & 0.918 &  0.908 & & 0.876 & 0.754         \\\cline{1-12}
\end{tabular}
\end{table}

\subsection{A test of  Hermite rank/normality}
\label{sec6.4}
In a concluding numerical example, we wish to illustrate that data blocking has impacts for inference under LRD beyond the resampling.
One basic application of data blocks is for testing the null hypothesis  that the Hermite rank is $m=1$ for a transformed LRD process $X_t=G(Z_t)$ against the alterative $m>1$.  This type of assessment has practical value in application. For example, analyses  in financial econometrics can involve LRD models with assumptions about $m$  (cf.~\cite{BCL98}).
More generally, inference from sample averages under LRD may use normal theory
only if $m=1$ \cite{T75,T77}.  Even considering resampling approximations under LRD, the block bootstrap (i.e., full  data re-creation) becomes valid when $m=1$ \cite{L93}, while subsampling (i.e., small scale  re-creation) should be used instead if $m>1$ \cite{BT17,BW18,HHJ95}.

Based on data $X_1,\ldots,X_n$ from a LRD process $X_t=G(Z_t)$, a simple assessment of $H_0:m=1$ can be based on data blocks of $\ell$ as follows. The idea is to make averages (say) $W_i \equiv \sum_{j= 1}^{\ell} X_{j + \ell(i-1)}/\ell $  of  length $\ell$ blocks, $i=1,\ldots,b\equiv \lceil n/\ell\rceil$, as  in Section~\ref{sec2.2}, and
then check their agreement to normality.  Letting  $\Phi(\cdot)$ denote a standard normal cdf,  we compare the collection of residuals $R_i\equiv \Phi( (W_i-\bar{W} )/S_W)$  to a uniform$(0,1)$ distribution, where  $\bar{W},S_W$ are the average and standard deviation of  $\{W_1,\ldots,W_b\}$.
In a usual fashion, we can assess uniformity by applying a Kolmogorov–Smirnov statistic or an Anderson-Darling \cite{AD52}   statistic (e.g., $A \equiv -b - b^{-1}\sum_{i=1}^b (2i-1)\log \{R_{(i)}[1-R_{(b+1-i)}]\}$ for ordered $R_{(i)}$).
  Of course, the distribution of such a test statistic requires calibration under the null $H_0:m=1$.  However, a central limit theorem \cite{T79} for LRD
processes $X_t = G(Z_t)$ when $m=1$ gives
\begin{eqnarray}
\label{eqn:limitresult}   \frac{1}{\sqrt{\var(\bar{X}_n)}} \frac{1}{n}\sum_{i = 1}^{\lfloor nt \rfloor}\left( X_i - \mu \right) \overset{d}{\to} c B_H(t), \quad 0\leq t\leq 1,
\end{eqnarray}
as $n\to\infty$, where $B_H(t)$ denotes   fractional Brownian motion with Hurst index $H= 1 - \alpha/2\in (0,1)$ and $c>0$ is a process constant.  Note that (\ref{eqn:limitresult})
no longer holds under LRD when $m>1$, which aids in  testing.  The property  (\ref{eqn:limitresult})
 suggests a simple bootstrap procedure for re-creating the null distribution of residual-based test statistics, given in Algorithm~\ref{tab:algo},   because such statistics are invariant to the  location-scale used in a bootstrap sample.

\begin{algorithm}[t]
\SetAlgoLined
\KwData{Given a LRD sample $X_1, .., X_n$.}
 Set initializations: block size $\ell$; number $M$ of resamples; and significance level $\alpha_{\mathrm{sig}} = 0.05$.

 Step 1. Calculate  a test statistic $T_0$ for normality (e.g., Anderson-Darling) from block averages\;

 Step 2. Estimate $\widehat{\alpha m}$ the memory parameter $\alpha m$ or Hurst Index $\widehat{H} = 1-\widehat{\alpha m}/2$\;

 Step 3. \For{$k =1,\ldots,M$}{
  (i) Simulate a Fractional Brownian motion sample $\{B^{*}_{\scriptscriptstyle \widehat{H}}(\frac{1}{n}),\ldots, B^{*}_{\scriptscriptstyle \widehat{H}}(\frac{n}{n})$\} with Hurst index $\widehat{H}$\;

  (ii) Obtain a bootstrap sample $X_1^*, \ldots, X_n^*$ as $X_1^{*} = B^{*}_{\scriptscriptstyle \widehat{H}}(\frac{1}{n})$ and $X_j^{*} = B^{*}_{\scriptscriptstyle \widehat{H}}(\frac{j}{n})- B^{*}_{\scriptscriptstyle \widehat{H}}(\frac{j-1}{n})$, $j = 2,\ldots, n$\;

  (iii) Calculate the $k$th bootstrap test statistic, $T^{*}_k$, for normality  from block averages in $X_1^*,\ldots,X_n^*$\;
 }

 Step 4. Compute $\widehat{q}$ as the $1-\alpha_{\mathrm{sig}}$ sample percentile of $\{T^{*}_1,\ldots,T^*_M\}$\;

 Step 5. Reject if $T_0 > \widehat{q}$.
 \caption{Hermite rank test for $m=1$ (normality) from  LRD series}
 \label{tab:algo}
\end{algorithm}

 The role of data blocking for tests of Hermite rank $H_0: m = 1$  with LRD series $X_t=G(Z_t)$ may be
 traced to recent work of \cite{BMG16}.  Those authors test for $m=1$ (normality) with a cumulant-based two-sample $t$-test, using two   samples  generated from data by different OL block resampling approaches.  Our block-based test  is different and perhaps more basic.
To briefly compare these tests,  we use data generation settings  from   \cite{BMG16} where $X_t = G(Z_t)$ with $G(z)= z+ 20^{-1}H_2(z)+ (20\sqrt{3})^{-1}H_3(z) $ (i.e., $m=1$) or $G(z)= \cos(z)$  (i.e., $m=2$), and $Z_t$ denotes a standardized FAIRMA$(0,(1-\alpha)/2,0)$ Gaussian process for $\alpha=0.2,0.8$.  The test in \cite{BMG16}
uses block lengths $\ell=n^{1/4}$ or $\ell=n^{1/2}$, where some best-case results  provided there
 assume the memory parameter to be known. To facilitate comparison against these, we simply use a similar block  $\ell=n^{1/2}$ for our test, as a reasonable  choice under LRD, and consider both OL/NOL blocks; we also use local Whittle estimation of the memory parameter along with 200 resamples in Algorithm~\ref{tab:algo}.
Table~\ref{tab:htestad} lists power (based on 500 simulation runs) of our test using a 5\% nominal level  compared to test findings of \cite{BMG16} (Table 2);
   we report an Anderson-Darling statistic in Table~\ref{tab:htestad} though a Kolmogorov–Smirnov statistic  produced similar results.
Both the proposed test and  the \cite{BMG16}-test maintain the nominal size for the LRD processes with $m=1$, but our block-based
test has much larger power for the LRD process defined by $m=2$ and $\alpha=0.2$.  The process defined by $m=2$ and $\alpha=0.8$ in
Table~\ref{tab:htestad} is actually SRD; as both our test and  the \cite{BMG16}-test are block-based assessments of normality, both
tests should maintain their sizes in this case and our test performs a bit better.
This illustrates that data-blocking has potential for assessments beyond usage in resampling.

\begin{table}[t]
\centering
\caption{Rejection rates of nominal 5\% tests of Hermite rank $m=1$   from block-based  statistic (Algorithm~\ref{tab:algo}
with OL/NOL blocks   $\ell = \lfloor n^{1/2} \rfloor$)  with LRD series $X_t\equiv G(Z_t)$ using
$G_1(z) = z+ 20^{-1}H_2(z)+ (20\sqrt{3})^{-1}H_3(z)$   or $G_2(z) = \cos(z)$. Results from the test of \cite{BMG16} (OL blocks $\ell = n^{1/2} $ or $ n^{1/4} $) are included.}

\label{tab:htestad}
\begin{tabular}{llrllccccccc}
\\
\cline{1-12}
               &   & &  &    & \multicolumn{3}{c}{$\alpha =0.20$} & & \multicolumn{3}{c}{$\alpha =0.80$}                                                       \\ \cline{6-12}
$G/m$            &    &\multicolumn{2}{c}{Testing Method}  & & $n=400$ & $ 1000$ & $ 10000$                                 & & $n=400$ &  1000 &  10000                               \\ \cline{1-12}
 \multirow{4}*{\parbox{1cm}{$G_1$,\\$m=1$}} & &\multirow{2}{*}{Algorithm~\ref{tab:algo}}& OL  & & 0.06  & 0.06  & 0.07                                   & & 0.05  &  0.08  & 0.07                                \\
                    & & &  NOL & &    0.10  &  0.06  & 0.08               &   &  0.06  & 0.05  & 0.06                                     \\\cline{3-12}

  & & \multirow{2}{*}{\cite{BMG16}-test}& $\ell =   n^{1/2} $   & & 0.04 & 0.05 & 0.03              & & 0.01 & 0.03 &  0.03                               \\
                   & &   & $\ell =  n^{1/4} $  & &   0.03 & 0.04 & 0.04             &   & 0.07 & 0.06 & 0.03                                  \\\hline\hline

 \multirow{4}*{\parbox{1cm}{$G_2$,\\$m=2$}} & & \multirow{2}{*}{Algorithm~\ref{tab:algo}}&  OL  & & 0.60 & 0.85 & 0.99                & &  0.13  & 0.14  & 0.13                               \\
                   & & &  NOL & &   0.45 &  0.68  & 0.99               &   &  0.06 & 0.09  & 0.07                                  \\\cline{3-12}

  & & \multirow{2}{*}{\cite{BMG16}-test}& $\ell =   n^{1/2} $  & & 0.11 & 0.23 & 0.26              & &  0.13  & 0.13 & 0.20                              \\
                   & & & $\ell =   n^{1/4} $ & &   0.19 &  0.24  & 0.32               &   &  0.24 & 0.16  & 0.27                                  \\\cline{1-12}
\end{tabular}
\end{table}



\section{Concluding Remarks}
\label{sec:7}
While block-based resampling methods  provide useful estimation under data dependence, their performance
is intricately linked to a  block length parameter, which is important  to understand.
This problem has been extensively investigated under weak or short-range dependence (SRD) (cf.~\cite{L03}, ch.~3), though relatively little has
been known for long-range dependence (LRD), especially outside  the pure Gaussian case (cf.~\cite{KN11}).
For general long-range dependent (LRD)   $X_t=G(Z_t)$ processes, based on  transforming a LRD Gaussian series $Z_t$,
results here showed that  properties and best block sizes with resampling variance estimators under LRD can intricately depend on
 covariance strength and the structure of non-linearity in  $G(\cdot)$.  The long-memory guess $O(n^{1/2})$ for   block size   \cite{HJL98} may have optimal order at times,  owing more to such non-linearity  rather than intuition  about LRD.  Additionally, we provided a data-based rule for best block selection,
 which was shown to be consistent under complex cases for blocks with LRD. While we   focused on a variance estimation problem with resampling under LRD, block selection for
  distribution estimation is also of interest, though   seemingly requires further  difficult study of distributional expansions for  statistics from LRD series $X_t=G(Z_t)$.  However, we showed that the block selections developed here can provide helpful benchmarks for choosing block size with resampling or other block-based inference problems under LRD.
  
  The current work may also suggest future
possibilities toward estimating the Hermite $m$ rank (or other ranks) under LRD.  No estimators of $m$ currently exist; instead,
only estimation of the long-memory exponent $\am$ of $X_t = G(Z_t)$ has been possible, which depends on the covariance decay rate $\alpha<1/m$ of $Z_t$.
Results here established that the variance of a block resampling estimator depends only on $\alpha$, apart from the Hermite rank $m$ itself.
  This suggests that some higher-order moment estimation may be investigated for separately estimating  the memory coefficient $\alpha$ and Hermite rank $m$
   under LRD.

\begin{appendix}
\section{Coefficient of Optimal Block Size}\label{appn}

The coefficient $K_{\alpha, m, m_{2}}$ of Corollary~\ref{cor1} is presented in cases   with notation:   $A  \equiv  \mathcal{B}_0^2(m)$, $B  \equiv (2  C_0^{m_2} J_{m_2}^2 /m_2!  )^2 $, $C \equiv  (\mathcal{B}_0(m) + 2  \sum_{j=m_2}^\infty     \sum_{k=1}^\infty   [\gamma_Z(k)]^j  J^2_j/j! )^2$,
  and $D  \equiv   ( 2  C_0^{m_2}/[(1-m_2 \alpha)(2-m_2 \alpha)] (J_{m_2}^2/m_2!))^2$ related to
  (\ref{eqn:B1});    $E \equiv v^2_{\infty, \alpha m}$ in (\ref{xbar.var}); and $F\equiv a_\alpha$ from Theorem~\ref{theorem2}.

\begin{description}
    \item Case $1$: $m_2 = \infty$
\begin{equation*}
K_{\alpha, m, m_{2}}= \left\{
\begin{array}{ll}
\frac{-(1-2\alpha)\sqrt{A E} + \sqrt{(1-2\alpha)^2 {A E} + 4 \alpha(1 - \alpha) A(E +F)}}{2 \alpha (E + F)} & \text { if } 0<\alpha<0.5, m = 1 \\
\left(\frac{A}{F} \right)^{0.5} & \text { if } \alpha=0.5, m=1 \\
\left(\frac{2 A (1 - \alpha)}{F}\right)^{\frac{1}{3 - 2\alpha}} & \text { if } 0.5<\alpha<1, m=1 \\
\left(\frac{A(1 - \alpha m)}{F \alpha} \right)^{\frac{1}{2(1 + \alpha - \alpha m)}} & \text { if } 0<\alpha<\frac{1}{m}, m \ge 2
\end{array}\right.
\end{equation*}

\item Case $2$: $m_2 < \infty$ with $\alpha m_2 > 1$
\begin{equation*}
K_{\alpha, m, m_{2}}= \left\{
\begin{array}{ll}
\frac{- (1 - 2\alpha)\sqrt{C E} + \sqrt{(1-2\alpha)^2 C E + 4 \alpha(1 - \alpha) C(E + F)}}{2 \alpha (E+F)} & \text { if } \frac{1}{m_2}<\alpha<0.5, m = 1 \\
\left(\frac{C}{F} \right)^{0.5} & \text { if } \alpha=0.5, m=1 \\
\left(\frac{2 C (1 - \alpha)}{F} \right)^{\frac{1}{3 - 2\alpha}} & \text { if } \max\{\frac{1}{m_2}, 0.5\}<\alpha<1, m=1 \\
\left(\frac{C(1 - \alpha m)}{F \alpha} \right)^{\frac{1}{2(1 + \alpha - \alpha m)}} & \text { if } \frac{1}{m_2}<\alpha<\frac{1}{m}, m \ge 2
\end{array}\right.
\end{equation*}

\item Case $3$: $m_2 < \infty$ with $\alpha m_2 = 1$

\begin{equation*}
K_{\alpha, m, m_{2}}= \left\{
\begin{array}{ll}
 \frac{- (1 - 2\alpha)\sqrt{B E} + \sqrt{(1-2\alpha)^2 B E + 4 \alpha(1 - \alpha) B(E + F)}}{2 \alpha (E+F)} & \text { if } 0<\alpha=\frac{1}{m_2}<0.5, m = 1 \\
\left(\frac{B}{F} \right)^{0.5} & \text { if } \alpha=\frac{1}{m_2} = 0.5, m=1 \\
\left(\frac{B(1 - \alpha m)}{F \alpha} \right)^{\frac{1}{2(1 + \alpha - \alpha m)}} & \text { if } 0<\alpha=\frac{1}{m_2}<0.5, m \ge 2
\end{array}\right.
\end{equation*}

\item Case $4$: $m_2 < \infty$ with $\alpha m_2 < 1$
\begin{equation*}
K_{\alpha, m, m_{2}}= \left\{
\begin{array}{ll}
\left(\frac{(m_2 -2) + \sqrt{(m_2 -2)^2 + (m_2-1)(E + F)D}}{E+F} \right)^{\frac{1}{\alpha m_2}} & \text { if } 0<\alpha<\frac{1}{m_2}, m = 1 \\
\left(\frac{D(m_2 - m)}{F} \right)^{\frac{1}{2 \alpha(1+m_2 -m)}} & \text { if } 0<\alpha<\frac{1}{m_2}, m \ge 2 \\
\end{array}\right.
\end{equation*}

\end{description}

\section{Proof of Theorem~\ref{theorem2}}\label{appn2}
The appendix considers the proof for the large-sample variance   (Theorem~\ref{theorem2})  of resampling estimators.  
Proofs are other results  are shown in the supplement \cite{ZLN22}.
To derive variance expansions of Theorem~\ref{theorem2},
we first consider the OL block variance estimator $\widehat{V}_{\ell,\am,\OL}= \sum_{i=1}^{N}\ell^\am(\bar{X}_{i,\ell}- \mu)^2/N -\ell^\am (\widehat{\mu}_{n,\OL}-\mu)^2  $ from (\ref{boot.var}), expressed in terms
of the process mean $\E X_t=\mu$, the number $N=n-\ell+1$ of blocks, the block averages $\bar{X}_{i,\ell} \equiv \sum_{j=i}^{i+\ell-1} X_j/\ell$ (integer $i$), and $\widehat{\mu}_{n,\OL}= \sum_{i=1}^{N}\bar{X}_{i,\ell}/N$.   Due to mean centering, we may assume
$\mu =0$  without loss of generality.
 We then write the variance of $\widehat{V}_{\ell,\am,\OL} $ as
\[
\var\left( \widehat{V}_{\ell,\am,\OL}  \right) = v_{1,\ell} + v_{2,\ell} - 2 c_{\ell},\qquad c_{\ell}  \equiv  \cov\left(\frac{1}{N}\sum_{i=1}^{N}\ell^\am \bar{X}_{i,\ell}^2, \ell^\am \widehat{\mu}_{n,\OL}^2    \right) = c_{a,\ell}+c_{b,\ell},
\]
\begin{eqnarray*}
v_{1,\ell} \equiv    \var\left( \frac{1}{N}\sum_{i=1}^{N}\ell^\am \bar{X}_{i,\ell}^2  \right)=v_{1a,\ell} + v_{1b,\ell},\qquad
v_{2,\ell}  \equiv   \var\left(\ell^\am \widehat{\mu}_{n,\OL}^2   \right) =  v_{2a,\ell} + v_{2b,\ell},
\end{eqnarray*}
where each variance/covariance component $v_{1,\ell}$, $v_{2,\ell}$ and $c_{\ell}$ is decomposed into two further subcomponents
\begin{eqnarray}\nonumber
v_{1a,\ell} &\equiv & 2 \frac{\ell^{2\am}}{N} \sum_{k=-N}^N \left(1-\frac{|k|}{N} \right) \left[\cov(\bar{X}_{0,\ell},\bar{X}_{k,\ell})\right]^2, \\
\label{eqn:vb} v_{1b,\ell} &\equiv & \frac{\ell^{2\am}}{N} \sum_{k=-N}^N \left(1-\frac{|k|}{N} \right) \mathrm{cum}(\bar{X}_{0,\ell}, \bar{X}_{0,\ell},\bar{X}_{k,\ell},\bar{X}_{k,\ell}),\\
 \nonumber v_{2a,\ell} &\equiv &2\ell^{2\am} \left[ \var(\widehat{\mu}_{n,\OL}) \right]^2,  \\
 \nonumber v_{2b,\ell} & \equiv & \ell^{2\am}\mathrm{cum}(\widehat{\mu}_{n,\OL},\widehat{\mu}_{n,\OL},\widehat{\mu}_{n,\OL},\widehat{\mu}_{n,\OL}),\\
 \nonumber c_{a,\ell} &\equiv& 2 \frac{\ell^{2\am}}{N}   \sum_{i=1 }^N  \left[\cov(\bar{X}_{i,\ell}, \widehat{\mu}_{n,\OL}) \right]^2, \\
 \nonumber c_{b,\ell} &\equiv& \frac{\ell^{2\am}}{N^3} \sum_{i=1}^{N}\sum_{j=1}^N\sum_{k=1}^N  \mathrm{cum}(\bar{X}_{i,\ell}, \bar{X}_{i,\ell},\bar{X}_{j,\ell},\bar{X}_{k,\ell}),
\end{eqnarray}
 consisting of sums involving 4th order cumulants ($v_{1b,\ell}, v_{2b,\ell},  c_{b,\ell}$) or sums involving covariances ($v_{1a,\ell}, v_{2a,\ell},  c_{a,\ell}$).
 The second decomposition step follows from the product theorem for cumulants (e.g., $\cov(Y_1Y_2,Y_3Y_4) =   \cov(Y_1,Y_3)\cov(Y_2,Y_4)+\cov(Y_1,Y_4)\cov(Y_2,Y_3)+\mathrm{cum}(Y_1,Y_2,Y_3,Y_4)$ for arbitrary random variables with $\E Y_i=0$ and $\E Y_i^4<\infty$).  Note that $\E X_t^4=\E [G(Z_0)]^4<\infty$ and $G\in \overline{\mathcal{G}}_4(1)$ imply these variance components
  exist  finitely for any $n,\ell$ (see (S.9)
   or Lemma~3
  of the supplement ~\cite{ZLN22}).
Collecting terms, we have
\[
\var\left( \widehat{V}_{\ell,\am,\OL}  \right) = \Gamma_{\ell,\OL}  + \Delta_{\ell,\OL}
\]
where $
\Gamma_{\ell,\OL}  \equiv v_{1a,\ell} + v_{2a,\ell}-2c_{a,\ell}$ and
$\Delta_{\ell,\OL}  \equiv v_{1b,\ell} + v_{2b,\ell}-2c_{b,\ell}$
denote sums over covariance-terms $\Gamma_{\ell,\OL}$ or sums of 4th order cumulant terms $\Delta_{\ell,\OL}$.
In the NOL estimator case $\widehat{V}_{\ell,\am,\NOL}$, the variance expansion is  similar
$
\var\left( \widehat{V}_{\ell,\am,\NOL}  \right) = \Gamma_{\ell,\NOL}  + \Delta_{\ell,\NOL}
$ with the convention that $\Gamma_{\ell,\NOL}, \Delta_{\ell,\NOL}$ are defined by replacing the OL block number $N=n-\ell+1$, averages
$\bar{X}_{i,\ell}$ (or $\bar{X}_{j,\ell},\bar{X}_{k,\ell}$) and estimator $\widehat{\mu}_{n,\OL}= \sum_{i=1}^{N}\bar{X}_{i,\ell}/N$ with the NOL counterparts $b=\lfloor n/\ell \rfloor$, $\bar{X}_{1+(i-1)\ell,\ell}$ (or $\bar{X}_{1+(j-1)\ell,\ell},\bar{X}_{1+(k-1)\ell,\ell}$) and
$\widehat{\mu}_{n,\NOL}= \sum_{i=1}^{b}\bar{X}_{1+(i-1)\ell,\ell}/b$ in $v_{1a,\ell},v_{2a,\ell},c_{a,\ell}, v_{1b,\ell}, v_{2b,\ell},c_{b,\ell}$.

  Let $\Gamma_{\ell}$ denote either counterpart  $\Gamma_{\ell,\OL}$ or $\Gamma_{\ell,\NOL}$, and $\Delta_\ell$ denote either   $\Delta_{\ell,\OL}$ or $\Delta_{\ell,\NOL}$.  Theorem~\ref{theorem2} then follows by establishing that
     \[
    \Gamma_{\ell} =\left\{   \begin{array}{lcl}
    O\left( \left(\frac{\ell}{n}\right)^{ \min\{ 1, 2\am\} } [\log n]^{I(2\am =1)}  \right)=  o\left(\left(\frac{\ell}{n}\right)^{2\alpha} \right) & & \mbox{if $m \geq 2$}\\[.1cm]
     a_\alpha \left(\frac{\ell}{n}\right)^{\min\{1,2\alpha\}} [\log n]^{I(\alpha =1/2)} \big(1+o(1)\big) && \mbox{if $m=1$},
                                 \end{array}
        \right.
   \]
and
  \[
    \Delta_{\ell} -  r_{n,\alpha,m,m_p} =\left\{   \begin{array}{lcl}
  \phi_{\alpha,m}  \left(\frac{\ell}{n}\right)^{2\alpha} \big(1+o(1)\big) & & \mbox{if $m \geq 2$}\\[.2cm]
       o \left( \left(\frac{\ell}{n}\right)^{\min\{1,2\alpha\}} [\log n]^{I(\alpha =1/2)} \right)&& \mbox{if $m=1$},
                                 \end{array}
        \right.
   \]
 where $I(\cdot)$  denotes an indicator function  and
  $r_{n,\alpha,m,m_p}$ is defined in Theorem~\ref{theorem2}. For reference, when the Hermite rank $m\geq 2$, the contribution of $\Delta_{\ell} \propto (\ell/n)^{2\alpha}$
dominates the variance of  $\widehat{V}_{\ell,\am,\OL}$ or $\widehat{V}_{\ell,\am,\NOL}$; when $m=1$, the contribution of $\Gamma_{\ell} \propto (\ell/n)^{\min\{1,2\alpha\}} [\log n]^{I(\alpha =1/2)}$ instead dominates the variance in Theorem~\ref{theorem2}.

To establish these expansions of  $\Gamma_{\ell},\Delta_{\ell}$, we require a series of technical lemmas (Lemmas~1-4),
involving certain graph-theoretic moment expansions.  To provide some illustration, Lemma~1
and its proof are outlined in Appendix~\ref{appn3}; the remaining lemmas are described in the supplement \cite{ZLN22}.  Define an order constant $\tau_{\ell,m} \equiv (\ell/n)^{2 \alpha}$ if $m \geq 2$ and $\tau_{\ell,m}\equiv (\ell/n)^{\min\{1,2\alpha\}} [\log n]^{I(\alpha =1/2)}$ if $m=1$; we suppress the dependence of $\tau_{\ell,m}$ on $n$ and $\alpha$ for simplicity.
Then, the above expansion of  $\Gamma_{\ell} $ follows directly from Lemma~4
(i.e., $\Gamma_{\ell} = a_\alpha \tau_{\ell,1} (1+o(1))$ if $m=1$ and $\Gamma_{\ell}=o(\tau_{\ell,m})$ if $m \geq 2$).
For handling $\Delta_{\ell}$, Lemma~1
gives that
$v_{1b,\ell} = r_{n,\alpha,m,m_p} + \phi_{\alpha,m}  \tau_{\ell,m} \big(1+o(1)\big)$ when $ m \geq 2$
and $v_{1b,\ell} = r_{n,\alpha,m,m_p} + o(\tau_{\ell,m}) $ when $m=1$.
Combined with this, the expansion of $\Delta_{\ell}$ then follows from Lemmas~2 and 3,
which respectively show that   $c_{b,\ell} = o(\tau_{\ell,m})$ and $v_{2b,\ell} = o(\tau_{\ell,m})$ for any $m \geq 1$.   $\Box$

\section{Lemma~\ref{lem1}  (dominant 4th order cumulant terms)}\label{appn3}
In the proof of Theorem~\ref{theorem2} (Appendix~\ref{appn2}), recall $v_{1b,\ell}$ from (\ref{eqn:vb}) represents a sum of 4th order cumulants block averages
with OL blocks, where the version with  NOL blocks is
   $v_{1b,\ell} \equiv  b^{-1}\ell^{2\am} \sum_{k=-b}^b \left(1-\frac{|k|}{b} \right) \mathrm{cum}(\bar{X}_{0,\ell}, \bar{X}_{0,\ell},\bar{X}_{ k\ell,\ell},\bar{X}_{k\ell,\ell})$.  Lemma~\ref{lem1} provides an expansion of  $v_{1b,\ell}$ under LRD, which is valid in either OL/NOL block case.

\begin{Lemma}
\label{lem1} Suppose the assumptions of Theorem~\ref{theorem2} ($G$ has Hermite rank $m \geq 1$    and Hermite pair-rank $m \leq m_p \leq \infty$) with positive constants  $ \phi_{\alpha,m}, \lambda_{\alpha,m_p}$
  there. Then,
  \[
  v_{1b,\ell} =
   r_{n,\alpha,m,m_p} +
   \left\{   \begin{array}{lcl}
  \phi_{\alpha,m}  \left(\frac{\ell}{n}\right)^{2\alpha} \big(1+o(1)\big) & & \mbox{if $m \geq 2$}\\[.2cm]
       o \left( \left(\frac{\ell}{n}\right)^{\min\{1,2\alpha\}} [\log n]^{I(\alpha =1/2)} \right)&& \mbox{if $m=1$},
                                 \end{array}
        \right.
   \]
where $I(\cdot)$  denotes an indicator function  and  $r_{n,\alpha,m,m_p}$ is from Theorem~\ref{theorem2}.
  \end{Lemma}

\noindent {\sc Remark 2:}  The proof of Lemma~\ref{lem1} involves a standard, but technical, graph-theoretic representation of the 4th order cumulant among Hermite polynomials $H_{k_1}(Y_1),H_{k_2}(Y_2),H_{k_3}(Y_3)$, $H_{k_4}(Y_4)$  ($k_1,k_2,k_3,k_4 \geq 1$)    at a generic sequence  $(Y_1,Y_2,Y_3,Y_4)$ of  marginally  standard normal variables, with covariances $\E Y_i Y_j = r_{ij}=r_{ji}$ for $1 \leq i<j\leq 4$.  Namely, it holds that
\begin{equation}
\label{eqn:4cuma}
\mathrm{cum}[H_{k_1}(Y_1),H_{k_2}(Y_2),H_{k_3}(Y_3),H_{k_4}(Y_4)] = k_1! k_2! k_3! k_4! \sum_{A \in \mathcal{A}_c(k_1,k_2,k_3,k_4)} g(A)r(A),
\end{equation}
where above $\mathcal{A}_c(k_1,k_2,k_3,k_4)$ denotes the collection of all path-connected multigraphs from a generic set of four points/vertices
$p_1,p_2,p_3,p_4$, such that point $p_i$ has degree $k_i$ for $i=1,2,3,4$.  Each multigraph $A \equiv (v_{12},v_{13},v_{14},v_{23},v_{24},v_{34}) \in \mathcal{A}_c(k_1,k_2,k_3,k_4)$ is defined by distinct counts  $v_{ij}=v_{ji}\geq 0$, interpreted as the number of graph lines connecting points $p_i$ and $p_j$, $1 \leq i < j \leq 4$; no other lines are possible in $A$.   Then $g(A)  \equiv 1/[\prod_{1 \leq i<j \leq 4}  (v_{ij}!)]$   represents a so-called multiplicity factor,   while $r(A)\equiv   \prod_{1 \leq i<j \leq 4} r_{ij}^{v_{ij}}$ represents a weighted product of covariances among variables in $(Y_1,Y_2,Y_3,Y_4)$ (cf.\cite{T77}).  Membership $A   \in \mathcal{A}_c(k_1,k_2,k_3,k_4)$ requires
 degrees $k_i = \sum_{j: j \neq i } v_{ij}$ for $i=1,2,3,4$ (e.g., $k_2 = v_{12}+v_{23}+v_{24}$) as well as
 a path-connection in $A$ between  any  two points $p_i$ and $p_j$;   namely, for any given $1\leq i<j\leq 4$,
  an index sequence $i\equiv i_0,i_1,\ldots, i_m \equiv j \in\{1,2,3,4\}$ (for some $1 \leq m\leq 3$) must exist whereby $v_{i_{u-1} i_u}>0$ holds for each $u=1,\ldots,m$, entailing that consecutive points among $p_{i}=p_{i_0}, p_{i_1},\ldots, p_{i_m}=p_{j}$ are connected with lines in $A$.  In (\ref{eqn:4cuma}),
$\mathrm{cum}[H_{k_1}(Y_1),H_{k_2}(Y_2),H_{k_3}(Y_3),H_{k_4}(Y_4)] =0$ holds  whenever $\mathcal{A}_c(k_1,k_2,k_3,k_4)$ is empty;  given integers $k_1,k_2,k_3,k_4 \geq 1$,
 $\mathcal{A}_c(k_1,k_2,k_3,k_4)$ will be empty    if $k_1+k_2+k_3+k_4=2q$ fails to hold for some integer $q \geq 2$ with $\max_{i=1,2,3,4} k_i \leq q$.
See the supplement~\cite{ZLN22} for more background  and details on this graph-theoretic representation.\\

\noindent \textbf{Proof of Lemma~\ref{lem1}}.   We focus on the OL block version (\ref{eqn:vb}) of $v_{1b,\ell}$;  the NOL block case follows by the same
essential arguments, though the cumulant sums involved with NOL blocks are
less involved and simpler to handle.  We     assume $\mu=0$ for reference; the mean $\mu \equiv \E G(Z_t)$ does not impact the 4th order cumulants here.
 Define $\tau_{\ell,m} \equiv (\ell/n)^{2 \alpha}$ if $m \geq 2$ (i.e., $\alpha<1/2$) and $\tau_{\ell,m}\equiv (\ell/n)^{\min\{1,2\alpha\}} [\log n]^{I(\alpha =1/2)}$ if $m=1$ to describe the order of interest in Lemma~\ref{lem1} along cases $m\geq 2$ or $m=1$.

Using Lemma~3
(i.e., $\max_{0 \leq k \leq 2\ell}|\mathrm{cum}(\bar{X}_{0,\ell},\bar{X}_{0,\ell},\bar{X}_{k,\ell},\bar{X}_{k,\ell})|$ is $O(\ell^{-2\am}))$ for $m \geq 2$ or $O(\ell^{-\alpha-\min\{1,2\alpha\}}[\log \ell]^{I(\alpha=1/2)})$ for $m=1$), we may truncate   the sum $v_{1b,\ell}$ in (\ref{eqn:vb})
  as
  \begin{eqnarray}
  \label{eqn:cum1}
    v_{1b,\ell} =   v_{1b,\ell}^{trun} +  O\left( \frac{\ell^{1+I(m=1)(\alpha -\min\{1,2\alpha\})}}{n} [\log \ell]^{I(\alpha=1/2)} \right),
    \end{eqnarray}
     where the order term   is $o(\tau_{\ell,m})$.  Then, using that $G \in \overline{\mathcal{G}}_4(1)$ (cf.~(S.9)
     of the supplement~\cite{ZLN22}) along with the 4th order cumulant form (\ref{eqn:4cuma}) for Hermite polynomials, we may use the multi-linearity of cumulants and the Hermite expansion (\ref{eqn:herm})  to express $v_{1b,\ell}^{trun}$ in (\ref{eqn:cum1}) as
     \begin{eqnarray}
   \label{eqn:cum}  &\quad&  2 \frac{\ell^{2\am}}{N}  \sum_{k=2\ell +1 }^N \left(1-\frac{k}{N} \right) \mathrm{cum}(\bar{X}_{0,\ell},\bar{X}_{0,\ell},\bar{X}_{k,\ell},\bar{X}_{k,\ell}) \;\equiv \;v_{1b,\ell}^{trun}\\
   \nonumber \quad & =&
     2 \frac{\ell^{2\am}}{N}  \sum_{k=2\ell +1 }^N \left(1-\frac{k}{N} \right)
     \sum_{m \leq k_1,k_2,k_3,k_4} \left(\prod_{i=1}^4 \frac{J_{k_i}}{k_i!}\right) \frac{1}{\ell^4}\times\\
     \nonumber  && \qquad \sum_{1\leq i_1,i_2,i_3,i_4 \leq \ell}
      \mathrm{cum}(H_{k_1}(Z_{i_1}),H_{k_2}(Z_{i_2}),H_{k_3}(Z_{i_3+k}), H_{k_4}(Z_{i_4+k}))\\
     \nonumber \quad & = &
     2 \sum_{q =2m}^\infty \sum_{k_1+k_2+k_3+k_4=2q,\atop
     m \leq k_1,k_2,k_3,k_4 \leq q}   \left(\prod_{i=1}^4  J_{k_i} \right) \sum_{A \in \mathcal{A}_c(k_1,k_2,k_3,k_4),\atop \mathcal{A}_c(k_1,k_2,k_3,k_4) \neq \emptyset} g(A)
     \frac{1}{N}\sum_{k=2\ell+1}^N \left(1-\frac{k}{N} \right) R_{\ell,k}(A)\\
   \nonumber \quad &\equiv &  t_{1,\ell} + t_{2,\ell} + t_{3,\ell} \mbox{\quad (say)},
     \end{eqnarray}
  with $g(A) =1/[\prod_{1 \leq i <j \leq 4} (v_{ij}!)]$, $A\equiv (v_{12},v_{13},v_{14},v_{23},v_{24},v_{34}) \in \mathcal{A}_c(k_1,k_2,k_3,k_4)$, and
     \begin{eqnarray}
     \label{eqn:RA}
   \qquad  R_{\ell,k}(A)
     &\equiv & \frac{\ell^{2 \am}}{\ell^4}\sum_{i_1=1}^{\ell}\sum_{i_2=1}^{\ell}\sum_{i_3=1}^{\ell}\sum_{i_4=1}^{\ell}
     [\gamma_Z(i_1-i_2)]^{v_{12}}  [\gamma_Z(i_1-i_3-k)]^{v_{13}}  \times\\
   \nonumber  &&  [\gamma_Z(i_1-i_4-k)]^{v_{14}} [\gamma_Z(i_2-i_3-k)]^{v_{23}} [\gamma_Z(i_2-i_4-k)]^{v_{24}} [\gamma_Z(i_3-i_4)]^{v_{34}}\end{eqnarray}
     in terms of Gaussian process covariances $\gamma_Z(\cdot)$ from (\ref{eqn:zcov}).  In (\ref{eqn:cum}), we use that the Hermite rank of $G$ is
     $m$ (so   $J_{k}=0$ for $k<m$); that the 4th order cumulants (\ref{eqn:4cuma}) of Hermite polynomials  are zero when a collection  $\mathcal{A}_c(k_1,k_2,k_3,k_4)$ is empty; and, relatedly, for given integers $k_1,k_2,k_3,k_4 \geq m$, a non-empty collection $\mathcal{A}_c(k_1,k_2,k_3,k_4)$   requires the number of lines, say   $q \equiv (k_1+k_2+k_3 + k_4)/2$, of a graph $A$ in $\mathcal{A}_c(k_1,k_2,k_3,k_4)$ to be an integer $q \geq 2m$ with $m \leq \max_{i=1,2,3,4}k_i \leq q$.
     The three components $t_{1,\ell},t_{2,\ell}$ and $t_{3,\ell}$ in (\ref{eqn:cum}) are defined by splitting the sum into three mutually exclusive cases, depending on the number of lines
       $q \equiv (k_1+k_2+k_3 + k_4)/2 = \sum_{1 \leq i < j\leq 4} v_{ij}$ and
       the value of $v_{13}+v_{14}+v_{23}+v_{24}=q-v_{12}-v_{34}$ for a  multigraph  $A\equiv (v_{12},v_{13},v_{14},v_{23},v_{24},v_{34})   \in \mathcal{A}_c(k_1,k_2,k_3,k_4) \neq \emptyset$; note that counts $v_{13},v_{14},v_{23},v_{24}$ involve covariances $\gamma_Z(\cdot)$
         at large lags in (\ref{eqn:RA}) (i.e., larger than $\ell$ by $2k \geq \ell+1$), which is not true of counts $v_{12},v_{34}$.  The three cases for defining $t_{1,\ell},t_{2,\ell}$ and $t_{3,\ell}$ are given by:
         (i) the case $q=2m$ with $k_1=k_2=k_3=k_4=m$, which yields
     \begin{equation}
     \label{eqn:t1start}
     t_{1,\ell} \equiv 2  J_m^4 \sum_{A \in \mathcal{A}_c(m,m,m,m) \neq \emptyset} g(A)
     \frac{1}{N}\sum_{k=2\ell+1}^N \left(1-\frac{k}{N} \right) R_{\ell,k}(A);
     \end{equation}

\begin{figure}
\centering
{\small \begin{tikzpicture}[node distance={30mm}, thick, main/.style = {draw, circle}]
\node[main] (1) {$p_1$};
\node[main] (2) [above right of=1] {$p_2$};
\node[main] (3) [below right of=1] {$p_3$};
\node[main] (4) [above right of=3] {$p_4$};
\draw (1) -- node[midway,  xshift=-0.6cm, pos=0.7] {$1$} (2);
\draw (1) -- node[midway,  yshift= -0.3cm, pos=0.3] {$1$} (4);
\draw (2) -- node[midway,  xshift= 0.3cm, pos=0.3] {$1$} (3);
\draw (3) -- node[midway,  xshift= 0.6cm, pos=0.3] {$1$} (4);
\end{tikzpicture} }
\caption{An example of multigraph $A  \equiv (1,0,1,1,0,1) \in \mathcal{A}_c(m,m,m,m)$ for $m=2$.}
\label{fig:multigraph}
\end{figure}
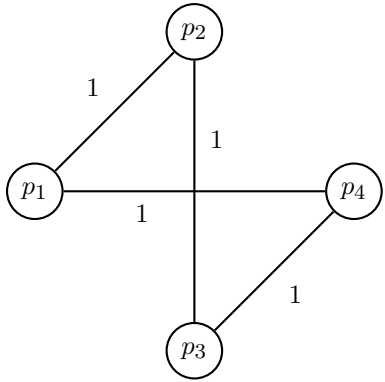

   (ii) the case that $q\geq 2m+1$ where the sum over $A \in \mathcal{A}_c(k_1,k_2,k_3,k_4)$ is also restricted to multigraphs $A\equiv (v_{12},v_{13},v_{14},v_{23},v_{24},v_{34})$ with $v_{13}+v_{14}+v_{23}+v_{24}=1$, which yields
     \[
    t_{2,\ell} \equiv  2 \sum_{q =2m+1}^\infty \sum_{k_1+k_2+k_3+k_4=2q,\atop
     m \leq k_1,k_2,k_3,k_4 \leq q}   \left(\prod_{i=1}^4  J_{k_i} \right) \sum_{A \in \mathcal{A}_c(k_1,k_2,k_3,k_4) \neq \emptyset,\atop
     v_{13}+v_{14}+v_{23}+v_{24}=1} g(A)
     \frac{1}{N}\sum_{k=2\ell+1}^N \left(1-\frac{k}{N} \right) R_{\ell,k}(A);
     \]
      and (iii) the final case that $q \geq 2m+1$ with the sum over $A \in \mathcal{A}_c(k_1,k_2,k_3,k_4)$ containing those connected graphs $A$ where $v_{13}+v_{14}+v_{23}+v_{24}\geq 2$, which yields
      \[
    t_{3,\ell} \equiv  2 \sum_{q =2m+1}^\infty \sum_{k_1+k_2+k_3+k_4=2q,\atop
     m \leq k_1,k_2,k_3,k_4 \leq q}   \left(\prod_{i=1}^4  J_{k_i} \right) \sum_{A \in \mathcal{A}_c(k_1,k_2,k_3,k_4) \neq \emptyset,\atop
     v_{13}+v_{14}+v_{23}+v_{24}\geq 2} g(A)
     \frac{1}{N}\sum_{k=2\ell+1}^N \left(1-\frac{k}{N} \right) R_{\ell,k}(A).
      \]
    Note that, for a connected multigraph $A$  here (i.e., $A\in  \mathcal{A}_c(k_1,k_2,k_3,k_4) \neq \emptyset$ for some $k_1,k_2,k_3,k_4\geq m$), a case $v_{13}+v_{14}+v_{23}+v_{24}=0$ is not possible (i.e.,  entailing that points $\{p_1,p_2\}$ are not path-connected in $A$ to points  $\{p_3,p_4\}$); consequently, terms $t_{2,\ell}$ and $t_{3,\ell}$ address all graphs $A$ with $q \geq 2m+1$ lines. Graphs with $q=2m$ lines appear in $t_{1,\ell}$.
     Lemma~\ref{lem1} will then follow from (\ref{eqn:cum1})-(\ref{eqn:cum}) by showing
     \begin{eqnarray}
     \label{eqn:t1}
      \qquad t_{3,\ell} &=&  o(\tau_{\ell,m}), \\\label{eqn:t2}
      \qquad  t_{2,\ell} & \equiv &  r_{n,\alpha,m,m_p}
   =I(m_p<\infty)\frac{\lambda_{\alpha,m_p}}{n^\alpha}
 \left(  \frac{\ell^{\am}}{\ell^{\min\{1,\alpha m_p\}}}  [\log \ell]^{I(m_p \alpha=1)}\right)^2
   \big(1+o(1)\big),\\
   \label{eqn:t3}
      \qquad t_{1,\ell} &=&  \phi_{\alpha,m}  \left(\frac{\ell}{n}\right)^{2\alpha} \big(1+o(1)\big) \quad\mbox{if $m \geq 2$} \qquad \& \qquad t_{1,\ell}=0 \quad\mbox{if $m=1$}.
    \end{eqnarray}
for $\tau_{\ell,m} \equiv (\ell/n)^{2 \alpha}$ if $m \geq 2$ and $\tau_{\ell,m}\equiv (\ell/n)^{\min\{1,2\alpha\}} [\log n]^{I(\alpha =1/2)}$ if $m=1$.

 We first consider showing (\ref{eqn:t2}) for $t_{2,\ell}$. Recall $t_{2,\ell}$ is defined by sums over connected multigraphs $A\equiv (v_{12},v_{13},v_{14},v_{23},v_{24},v_{34})\in \mathcal{A}_c(k_1,k_2,k_3,k_4) \neq \emptyset $ involving $q\equiv (k_1+k_2+k_3+k_4)/2 \geq 2m+1$ lines and $1=q-v_{12}-v_{34}=v_{13}+v_{14}+v_{23}+v_{24}$ (with $k_1,k_2,k_3,k_4 \geq m$).  In any such graph $A$, exactly one value among $v_{13},v_{14},v_{23},v_{24}$ equals 1, implying that four possible cases for the configuration of degrees, namely $(k_1 = v_{12},k_{4}=v_{34}, k_2=k_1+1, k_3=k_4+1)$ or $(k_1 = v_{12},k_{3}=v_{34}, k_2=k_1+1, k_4=k_3+1)$ or  $(k_2 = v_{12},k_{4}=v_{34}, k_1=k_2+1, k_3=k_4+1)$ or $(k_2 = v_{12},k_{3}=v_{34}, k_1=k_2+1, k_4=k_3+1)$ with $k_1,k_2,k_3,k_4 \geq m$.  Consequently, $t_{2,\ell}$ can be re-written using (\ref{eqn:RA}) as
\begin{equation}
   \label{eqn:t2l}
   t_{2,\ell} =2  \sum_{j_1=m_p}^\infty  \sum_{j_2=m_p}^\infty  J_{j_1}J_{j_1+1} J_{j_2 } J_{j_2+1}  \frac{1}{j_1! j_2!}
     \frac{1}{N}\sum_{k=2\ell+1}^N \left(1-\frac{k}{N} \right) R_{\ell,k,j_1,j_2}^*,
   \end{equation}
   due the Hermite pair-rank $m_p \equiv \inf\{k \geq 1: J_{k}J_{k+1} \neq 0\} \geq m$, for a covariance sum
\begin{eqnarray*}
 R_{\ell,k,j_1,j_2}^* &\equiv &  \frac{\ell^{2 \am}}{\ell^4}\sum_{i_1=1}^{\ell}\sum_{i_2=1}^{\ell}\sum_{i_3=1}^{\ell}\sum_{i_4=1}^{\ell}
     [\gamma_Z(i_1-i_2)]^{j_1} [\gamma_Z(i_3-i_4)]^{j_2} \times \\
     &&\qquad \left(   \gamma_Z(i_1-i_3-k)
   +  \gamma_Z(i_1-i_4-k)  +  \gamma_Z(i_2-i_3-k)  +  \gamma_Z(i_2-i_4-k) \right),
\end{eqnarray*}
$k \geq 2\ell +1$.   Note that $t_{2,\ell}=0$ if $m_p = +\infty$ (i.e., if $J_{k}J_{k+1}=0$ for all $k \geq 1$).  Hence, we may assume $m_p<\infty$ in the following to establish the form of $t_{2,\ell}$  in (\ref{eqn:t2}) along the possibilities $\alpha m_p <1$, $\alpha m_p=1$, or $\alpha m_p>1$.

   Due to the Gaussian covariance decay $\gamma_Z(k) \sim C_0 k^{-\alpha}$ as $k \to \infty$ (cf.~(\ref{eqn:zcov})), note that, for $k \geq 2$, $0 \leq i \leq \ell-1$ and $j_1,j_2 \geq m_p$, we can express   $ R_{\ell,k\ell+i,j_1,j_2}^*$ in (\ref{eqn:t2l}) as
   \begin{equation}
   \label{eqn:t2l2}
    R_{\ell,k\ell+i,j_1,j_2}^*  =    \frac{\ell^{2 \am}}{\ell^{\min\{1, j_1 \alpha\}+\min\{1, j_2 \alpha\}}} L_{\ell,j_1} L_{\ell,j_2} \times 4 \ell^{-\alpha} C_0 k^{-\alpha} +  R_{\ell,k\ell+i,j_1,j_2}^{**},
   \end{equation}
   through a covariance-type average
   \begin{equation}
   \label{eqn:L}
    L_{\ell,v}  \equiv   \frac{\ell^{\min\{1, v \alpha\}}}{\ell}\sum_{j=-\ell}^\ell\left( 1 - \frac{|j|}{\ell} \right) [\gamma_Z(j)]^v, \quad v=0,1,2\ldots,
   \end{equation}
 along with a  remainder  $R_{\ell,k\ell+i,j_1,j_2}^{**}$ satisfying
   \[
     |R_{\ell,k\ell+i,j_1,j_2}^{**} | \leq     C\frac{\ell^{2 \am} (\log \ell)^{I(j_1 = 1/\alpha) + I(j_2 = 1/\alpha)}}{\ell^{\min\{1, j_1 \alpha\}+\min\{1, j_2 \alpha\}}} \ell^{-\alpha}(k-1)^{-\alpha}4 C_0\left(\varepsilon_\ell   + 2(k-1)^{-1}\right)
   \]
    for $\varepsilon_\ell \equiv   \sup_{j \geq \ell} |\gamma_{Z}(j)/[C_0 j^{-\alpha}] -1 | =o(1)$  and for constants $C,C_1 \geq 1$ not depending on $\ell$, $k \geq 2$, $0 \leq i \leq \ell-1$ or $j_1,j_2 \geq m_p$; here the bound on the remainder follows from
   \begin{equation}
   \label{eqn:lcov}
     \sup_{ v \in \{0,1,2,\ldots\}} (\log \ell)^{-I(v=1/\alpha)} \ell^{-1+\min\{1,v  \alpha\}}\sum_{t=-\ell}^\ell |\gamma_Z(t)|^{v}  \leq C
 \end{equation}
 for some $C>0$  (i.e., by applying $|\gamma_Z(t)| \leq \min\{1,C_1 |t|^{-\alpha}\}$ for all $t \geq 1$)  along with $|\ell k + i + i_1 - i_2|\geq \ell(k-1)$ and $|\ell^{\alpha}|\ell k + i + i_1 - i_2|^{-\alpha} - k^{-\alpha}|\leq 2 |k-1|^{-1-\alpha}$
    (by Taylor expansion) for any $k \geq 2$, $0 \leq i, i_1,i_2 \leq \ell$.     Hence, if $m_p \alpha < 1$, we use (\ref{eqn:t2l})-(\ref{eqn:lcov}),   with $ \sum_{j=m_p}^\infty |J_{j}J_{j+1} | /j! <\infty$ by $G \in \overline{\mathcal{G}}_4(1)$ (cf.~(S.9)
    of \cite{ZLN22}), to write
     \begin{eqnarray*}
      t_{2,\ell} &=&2  \left( \frac{J_{m_p}J_{m_p+1}}{m_p! }\right)^2  L_{\ell,m_p}^2  \frac{\ell^{2 \am}}{\ell^{2 \alpha m_p}}
     \frac{\ell^{1-\alpha}}{N}\sum_{k=2}^{\lceil n/\ell \rceil} (1-N^{-1}\ell k )  4 C_0  k^{-\alpha} + \frac{\ell^{2 \am}}{\ell^{2 \alpha m_p}}  O\left(  \frac{\ell^{1-\alpha}}{N}\right)\\
      &&
     \quad   +   \left( \frac{\ell^{2 \am}(\log \ell)^{2}}{\ell^{ \alpha m_p  + \min\{1,\alpha(1+m_p)\}}} +
      \frac{\ell^{2 \am}}{\ell^{2 \alpha m_p}}O\left(\frac{\ell}{n}+\varepsilon_\ell\right)\right) O\left(  \frac{\ell^{1-\alpha}}{N}\sum_{k=2}^{\lceil n/\ell \rceil}    k^{-\alpha}  \right)\\
      &=&  8 C_0  \left( \frac{J_{m_p}J_{m_p+1}}{m_p!} \frac{2 C_0^{m_p}}{(1-\alpha m_p)(2-\alpha m_p)}\right)^2  \frac{1}{(1-\alpha)(2-\alpha)}
      \frac{\ell^{2 \am}}{\ell^{2 \alpha m_p}} \frac{1}{n^\alpha} (1+o(1)),
     \end{eqnarray*}
  upon applying, as $n\to \infty$,  that  $\varepsilon_\ell\rightarrow 0$, $\ell/n\rightarrow 0$, $N/n \rightarrow 1$, $\ell^{\alpha m_p - \min\{1,\alpha(1+m_p)\}} (\log \ell)^2\rightarrow 0$  for $m_p\alpha <1$ ($0<\alpha<1$) and
    \[
    L_{\ell,m_p}\rightarrow \frac{2 C_0^{m_p}}{(1-\alpha m_p)(2-\alpha m_p)},\quad \left(\frac{\ell}{N}\right)^{1-\alpha}\sum_{k=2}^{\lceil n/\ell \rceil} (1-N^{-1}\ell k )k^{-\alpha}\rightarrow \frac{1}{(1-\alpha)(2-\alpha)}.
    \]
    This shows (\ref{eqn:t2}) for $t_{2,\ell}$ when   $\alpha m_p <1$.
     If $\alpha m_p=1$ holds, then the derivation is similar by    (\ref{eqn:t2l})-(\ref{eqn:lcov}), where
      \begin{eqnarray*}
      t_{2,\ell} &=&2  \left( \frac{J_{m_p}J_{m_p+1}}{m_p! }\right)^2  \frac{L_{\ell,m_p}^2}{(\log \ell)^2}  \frac{\ell^{2 \am}}{\ell^{2 }} (\log \ell)^2
     \frac{\ell^{1-\alpha}}{N}\sum_{k=2}^{\lceil n/\ell \rceil} (1-N^{-1}\ell k )  4 C_0  k^{-\alpha} \\
      &&+\frac{\ell^{2 \am}(\log \ell)^2 }{\ell^{2 }} O\left(  \frac{\ell^{1-\alpha}}{N}\right)+
       \frac{\ell^{2 \am}(\log \ell)^2 }{\ell^{2 }} \left(1 +
     O\left(\frac{\ell}{n}+\varepsilon_\ell\right)\right) O\left(  \frac{\ell^{1-\alpha}}{N}\sum_{k=2}^{\lceil n/\ell \rceil}    k^{-\alpha}  \right)\\
      &=&  8 C_0  \left( \frac{J_{m_p}J_{m_p+1}}{m_p!} 2 C_0^{m_p} \right)^2\frac{1}{(1-\alpha)(2-\alpha)}
      \frac{\ell^{2 \am}}{\ell^{2  }}   (\log \ell)^2 \frac{1}{n^\alpha} (1+o(1)),
     \end{eqnarray*}
     using instead
      \[
      \frac{L_{\ell,m_p}}{\log \ell}\rightarrow  2 C_0^{m_p}>0;
     \]
  this then gives (\ref{eqn:t2}) for $t_{2,\ell}$ when   $\alpha m_p =1$.   In the final case $\alpha m_p>1$, we analogously have
    \begin{eqnarray*}
      t_{2,\ell} &=&2  \left(  \sum_{j=m_p}^\infty    \frac{J_{j}J_{j+1}}{j! } L_{\ell,j} \right)^2   \frac{\ell^{2 \am}}{\ell^{2}}
     \frac{\ell^{1-\alpha}}{N}\sum_{k=2}^{\lceil n/\ell \rceil} (1-N^{-1}\ell k )  4 C_0  k^{-\alpha}\\
      &&
       + \frac{\ell^{2 \am}}{\ell^{2 }}   \left(O\left(  \frac{\ell^{1-\alpha}}{N}\right)+
      O\left(\frac{\ell}{n}+\varepsilon_\ell\right)O\left(  \frac{\ell^{1-\alpha}}{N}\sum_{k=2}^{\lceil n/\ell \rceil}    k^{-\alpha}  \right)\right)    \\
      &=&  8 C_0  \left(  \sum_{j=m_p}^\infty    \frac{J_{j}J_{j+1}}{j! } \sum_{k=-\infty}^{\infty} [\gamma_Z(k)]^j \right)^2
      \frac{1}{(1-\alpha)(2-\alpha)} \frac{\ell^{2 \am}}{\ell^{2}}   \frac{1}{n^\alpha} (1+o(1)),
     \end{eqnarray*}
     where
     \[
     \sum_{j=m_p}^\infty    \frac{J_{j}J_{j+1}}{j! } L_{\ell,j} \rightarrow \sum_{j=m_p}^\infty    \frac{J_{j}J_{j+1}}{j! }  \sum_{k=-\infty}^{\infty} [\gamma_Z(k)]^j \in \mathbb{R}
     \]
     follows from the Dominated Convergence Theorem using that $L_{\ell,s}\rightarrow \sum_{k=-\infty}^{\infty} [\gamma_Z(k)]^v $ for each $v \geq m_p$; that $\sup_{v \geq m_p}|L_{\ell,v}| \leq
      \sum_{k=-\infty}^{\infty} |\gamma_Z(k)|^{m_p}<\infty$ holds for all $\ell$ (i.e., $\alpha m_p>1$ and $|\gamma_Z(k)|^{v} \leq |\gamma_Z(k)|^{m_p} \leq C_1^{m_p} k^{-\alpha m_p}$ for all $k\geq 1$, $v\geq m_p$ by $\gamma_Z(0)=1$); and that
         $ \sum_{j=m_p}^\infty |J_{j}J_{j+1} | /j! <\infty$ (by $G \in \overline{\mathcal{G}}_4(1)$).  This shows (\ref{eqn:t2}) for $t_{2,\ell}$ when   $\alpha m_p >1$ and concludes the establishment of (\ref{eqn:t2}).

 We next consider   showing (\ref{eqn:t1}) for $t_{3,\ell}$.  For Gaussian covariances $\gamma_Z(\cdot)$, using again that
 $|\gamma_Z(k)|\leq \min\{1, C_1 k^{-\alpha}\}$ for all $k \geq 1$ (for some $C_1 \geq 1$)   under (\ref{eqn:zcov}), we may bound
 \[
  \max_{0 \leq i,i_1,i_2 \leq \ell} |\gamma_Z(k\ell + i + i_1-i_2)| \leq \min\{1, C_1 \ell^{-\alpha} (k-1)^{-\alpha}\}
 \]
 for any $k \geq 2$, which follows by $|k\ell + i + i_1-i_2| \geq  \ell(k-1)$.  Applying this covariance inequality in $R_{\ell, \ell k+i}(A)$ from (\ref{eqn:RA}) shows that, for a multigraph $A\equiv (v_{12},v_{13},v_{14},v_{23},v_{24},v_{34})$ with $q=\sum_{1 \leq i < j\leq 4} v_{ij}$ lines,
 we have
 \begin{eqnarray}
 \nonumber
 && |R_{\ell, \ell k+i}(A)| \\ \nonumber &\leq& \frac{\ell^{2 \am}}{\ell^2} \left( \sum_{t=-\ell}^\ell |\gamma_Z(t)|^{v_{12}} \right)\left( \sum_{t=-\ell}^\ell |\gamma_Z(t)|^{v_{34}} \right) [\min\{1, C_1 \ell^{-\alpha} (k-1)^{-\alpha}\}]^{q-v_{12}-v_{34}}\\ \nonumber
 &\leq &   C  \frac{\ell^{2 \am}(\log \ell)^{I(v_{12}=1/\alpha)+I(v_{12}=1/\alpha)} }{\ell^{\min\{1,v_{12} \alpha\}+ \min\{1,v_{34} \alpha\}}}  [\min\{1, C_1 \ell^{-\alpha} (k-1)^{-\alpha}\}]^{q-v_{12}-v_{34}}\\
 \label{eqn:RA2} \qquad &\equiv &  B_{\ell,k,v_{12},v_{34},q-v_{12}-v_{34}}
 \end{eqnarray}
 holds for a generic constant $C>0$, not depending on the integers $k, \ell\geq 2$ or $0\leq i \leq \ell-1$, or the values of $v_{12},v_{34},q-v_{12}-v_{34}=v_{13}+v_{14}+v_{23}+v_{24} \in  \{0,1,2\ldots,\}$; above $I(\cdot)$ denotes an indicator function and we used  (\ref{eqn:lcov}) for establishing the bound in (\ref{eqn:RA2}).

      Now because
$G \in \overline{\mathcal{G}}_4(1)$ (cf.~(S.9)
of \cite{ZLN22}) and because
\[
 \frac{1}{N}\sum_{k=2\ell+1}^N \left(1-\frac{k}{N} \right) |R_{\ell,k}(A)| \leq \frac{\ell}{N} \sum_{k=2}^{\lceil n/\ell \rceil} B_{\ell,k,v_{12},v_{34},q-v_{12}-v_{34}}
\]
for a multigraph $A\equiv (v_{12},v_{13},v_{14},v_{23},v_{24},v_{34})$ with $q=\sum_{1 \leq i < j\leq 4} v_{ij}$ by (\ref{eqn:RA2}),  then $t_{3,\ell}=o(\tau_{\ell,m})$ will hold in (\ref{eqn:t1}) by establishing
\begin{equation}\label{eqn:B}
\sup_{A\equiv (v_{12},v_{13},v_{14},v_{23},v_{24},v_{34}): \atop q \geq 2m+1; q-v_{12}-v_{34} \geq 2;
k_1,k_2,k_3,k_4 \geq m} \frac{\ell}{n} \sum_{k=2}^{\lceil n/\ell \rceil} B_{\ell,k,v_{12},v_{34},q-v_{12}-v_{34}}= o(\tau_{\ell,m}),
\end{equation}
noting $n/N \rightarrow 1$ for $N=n-\ell+1$.
First, consider the case that $q- v_{12}-v_{34} \in \{2,\ldots,2m+1\}$ with $0\leq v_{12},v_{34}<1/\alpha$ so that  (\ref{eqn:RA2}) yields
\begin{eqnarray*}
 \frac{\ell}{n} \sum_{k=2}^{\lceil n/\ell \rceil}  B_{\ell,k,v_{12},v_{34},q-v_{12}-v_{34}} &\leq &C\frac{\ell^{2 \am}} {\ell^{ v_{12} \alpha  + v_{34}\alpha}}   \frac{\ell}{n} \sum_{k=1}^{\lceil n/\ell \rceil} C_1^{2m+1} (\ell k)^{-\alpha(q-v_{12}-v_{34})}\\
&\leq &C C_1^{2m+1} \frac{\ell^{2 \am-\alpha(q-v_{12}-v_{34})}} {\ell^{ v_{12} \alpha  + v_{34}\alpha}}  \frac{\ell}{n} \sum_{k=1}^{\lceil n/\ell \rceil} k^{-\alpha 2} =O(\ell^{-\alpha}) O(\tau_{\ell,m}),
\end{eqnarray*}
using above that $q \geq 2m+1$; that $k^{-\alpha(q-v_{12}-v_{34})} \leq k^{-\alpha 2}$ for $k\geq 1$; and that $\ell/n \sum_{k=1}^{\lceil n/\ell \rceil} k^{-\alpha 2}=O(\tau_{\ell,m})$ where $\ell^{-1}+\ell/n\rightarrow 0$. Note that the bound above is consequently $o(\tau_{\ell,m})$ and does not depend on the exact values of
$0\leq v_{12},v_{34}<1/\alpha$.
In the case that   $q- v_{12}-v_{34} \in \{2,\ldots,2m+1\}$ with $v_{12},v_{34} \geq  \lceil 1/\alpha \rceil$, we similarly obtain from (\ref{eqn:RA2}) that
\begin{eqnarray*}
 \frac{\ell}{n} \sum_{k=2}^{\lceil n/\ell \rceil}  B_{\ell,k,v_{12},v_{34},q-v_{12}-v_{34}} &\leq &C\frac{(\log \ell)^2 \ell^{2 \am} }{\ell^{2 }} \frac{\ell}{n}\sum_{k=1}^{\lceil n/\ell \rceil} C_1^{2m+1} (\ell k)^{-\alpha(q-v_{12}-v_{34})}\\
&\leq &C C_1^{2m+1} \frac{(\log \ell)^2 \ell^{2 \am-2\alpha}  }{\ell^{2 }}    \frac{\ell}{n} \sum_{k=1}^{\lceil n/\ell \rceil} k^{-\alpha 2}= O(\ell^{-\alpha}) O(\tau_{\ell,m}),
\end{eqnarray*}
using above that $ (\ell k)^{-\alpha(q-v_{12}-v_{34})} \leq (\ell k)^{-\alpha 2}$ for $k\geq 1$ and $\am \leq 1$; again the bound above is $o(\tau_{\ell,m})$ and does not depend on the exact values of  $v_{12},v_{34} \geq  \lceil 1/\alpha \rceil$.  When  $q- v_{12}-v_{34} \in \{2,\ldots,2m+1\}$ with $v_{12} \geq  \lceil 1/\alpha \rceil$ and $v_{34}<1/\alpha$, we use $q-v_{12}-v_{34}= v_{13}+v_{14}+v_{23}+v_{24}=k_3 + k_4 - 2 v_{34} \geq 2m -2 v_{34}$ in (\ref{eqn:B}) (i.e., $k_3,k_4 \geq m$) to write from (\ref{eqn:RA2}) that
 \begin{eqnarray*}
  \frac{\ell}{n} \sum_{k=2}^{\lceil n/\ell \rceil}  B_{\ell,k,v_{12},v_{34},q-v_{12}-v_{34}} &\leq &
  C\frac{(\log \ell)\ell^{2 \am}}{\ell^{1 + v_{34} \alpha}} \frac{\ell}{n} \sum_{k=1}^{\lceil n/\ell \rceil} C_1^{2m+1} (\ell k)^{-\alpha(q-v_{12}-v_{34})}\\
   &\leq &
   C C_1^{2m+1} \frac{(\log \ell) \ell^{2 \am} }{\ell^{1 + v_{34} \alpha}} \ell^{-\alpha(2m -2v_{34})}  \frac{\ell}{n} \sum_{k=1}^{\lceil n/\ell \rceil} k^{-2\alpha} \\&=&O( \ell^{ \alpha v_{*}-1} \log \ell ) O(\tau_{\ell,m}) =o(\tau_{\ell,m}),
\end{eqnarray*}
using $k^{-\alpha(q-v_{12}-v_{34})} \leq k^{-\alpha 2}$ for $k\geq 1$ as well as $\alpha v_{34}\leq \alpha v_{*}<1$ for $v_*= \lceil 1/\alpha \rceil -1$ as the largest integer less than $1/\alpha$;  the bound above does not depend on the exact values of  $v_{12} \geq  \lceil 1/\alpha \rceil$ and $v_{34}<1/\alpha$ and  the analog of the above also holds
similarly when $v_{34} \geq  \lceil 1/\alpha \rceil$ and $v_{12}<1/\alpha$ (using $q-v_{12}-v_{34}=  k_1 + k_2 - 2 v_{12} \geq 2m -2 v_{12}$).
Hence, (\ref{eqn:B}) will now follow by treating a final case that $q-v_{12}-v_{34} \geq 2m+2$.  When $q-v_{12}-v_{34} \geq 2m+2$, we use the bound in (\ref{eqn:RA2}) to write (for any $v_{12},v_{34} \geq 0$ and $q-v_{12}-v_{34} \geq 2m+2$) that
\begin{eqnarray*}
  \frac{\ell}{n} \sum_{k=2}^{\lceil n/\ell \rceil}  B_{\ell,k,v_{12},v_{34},q-v_{12}-v_{34}} &\leq &
  C\ell^{2\am} (\log \ell)^2 \frac{\ell}{n}  \sum_{k=1}^{\lceil n/\ell \rceil}   [C_1 (\ell k)^{-\alpha}]^{2m+2}= o(\tau_{\ell,m}).
\end{eqnarray*}
This establishes (\ref{eqn:B}) and consequently $t_{\ell,3}=o(\tau_{\ell,m})$ in (\ref{eqn:t1}).

         To complete the proof of Lemma~\ref{lem1}, we now consider establishing (\ref{eqn:t3}) for the term $t_{1,\ell}$, shown in
        (\ref{eqn:t1start}) involving a sum over
           multigraphs $A \equiv (v_{12},v_{13},v_{14},v_{23},v_{24},v_{34})\in \mathcal{A}_c(m,m,m,m)$.  For such $A$, the same degree requirement $k_1=k_2=k_3=k_4=m$ entails that $v_{12}=v_{34}$, $v_{13}=v_{24}$, $v_{14}=v_{23}$ so that we may prescribe the multigraph $A \equiv  (v_{12},v_{13},m-v_{12}-v_{13},m-v_{12}-v_{13},v_{13},v_{12})$ having $q=2m$ lines in terms of two counts $0\leq v_{12}, v_{13}$ where $v_{12}+v_{13}\leq m$.  Additionally, for a connected multigraph $A \in \mathcal{A}_c(m,m,m,m)$, it cannot be the case that $v_{12}=m$ holds, which would imply $v_{12}=m$ lines between points $\{p_1,p_2\}$ and $v_{34}=m$ lines between points $\{p_3,p_4\}$, with no lines between these two point groups (i.e., $A$ would not be connected); likewise, it cannot be the case that $v_{13}=m$ or that $v_{12}+v_{13}=0$ (so that $v_{14}=k_1-v_{12}-v_{13}=m$).  For this reason, $\mathcal{A}_c(m,m,m,m)$ is empty when $m=1$, so that the sum $t_{1,\ell}=0$  in (\ref{eqn:t3}) for $m=1$.  When $ m \geq 2$, it holds that $\mathcal{A}_c(m,m,m,m) \neq \emptyset$  and that the largest possible value of $v_{12}$ for some $A \in \mathcal{A}_c(m,m,m,m)$ is $m-1$, occurring when
     $v_{12}=m-1=v_{34}$ with either $v_{13}=1=v_{24}$ or $v_{14}=1=v_{23}$ (while, for reference, the smallest possible value of $v_{12}$ is $0$, occurring when $v_{12}=0=v_{34}$ and $1 \leq v_{13}=v_{24} <m$ with $v_{14}=m-v_{13}=v_{23}$).   As a consequence, when $m \geq 2$, we will split the sum $t_{1,\ell} = s_{1,\ell} + s_{2,\ell}$ into two parts, involving either a sum over $A \in \mathcal{A}_c(m,m,m,m)$ where $v_{12}=v_{34}= m-1 $  (given by $s_{1,\ell}$) or a sum over $A \in \mathcal{A}_c(m,m,m,m)$ where $0 \leq v_{12} \leq m-2 $ (given by $s_{2,\ell}$).

      For
     the first sum, we use the form of $g(A)\equiv 1/[\prod_{1 \leq i <j \leq 4}(v_{ij}!)]$ and $R_{\ell,k}(A)$ in $t_{1,\ell}$ (cf.~(\ref{eqn:RA})) to write
     \[
     s_{1,\ell} = 2  J_m^4 \left(\frac{1}{(m-1)!}\right)^2  N^{-1}\sum_{k=2\ell+1}^N  (1- N^{-1}k) R_{\ell,k}^\dag
     \]
    for
     \begin{eqnarray*}
     R_{\ell,k}^\dag  &\equiv &  \frac{\ell^{2 \am}}{\ell^4}\sum_{i_1=1}^{\ell}\sum_{i_2=1}^{\ell}\sum_{i_3=1}^{\ell}\sum_{i_4=1}^{\ell}
     [\gamma_Z(i_1-i_2) \gamma_Z(i_3-i_4)]^{m-1} \times\\
     &&\qquad [  \gamma_Z(i_1-i_3-k)\gamma_Z(i_2-i_4-k) +   \gamma_Z(i_1-i_4-k)\gamma_Z(i_2-i_3-k)], \quad k \geq  2\ell.
       \end{eqnarray*}
     Similarly, to the expansion in (\ref{eqn:t2l2}), we may use
  the Gaussian covariance decay $\gamma_Z(k) \sim C_0 k^{-\alpha}$ as $k \to \infty$ (cf.~(\ref{eqn:zcov})) and Taylor expansion to re-write $R_{\ell,k\ell+i}^\dag$   for any $k \geq 2$ and $0 \leq i \leq \ell-1$ as
   \begin{eqnarray*}
    R_{\ell,k\ell+i}^\dag& = &  \left(\frac{\ell^{\am}}{\ell^{\alpha (m-1)}} L_{\ell,m-1} \right)^2 2 \ell^{-2\alpha} C_0^{2} k^{-2\alpha} +  R_{\ell,k\ell+i}^{\ddag},
   \end{eqnarray*}
   with a covariance sum $L_{\ell,m-1}$ as in (\ref{eqn:L}) and a  remainder  $R_{\ell,k\ell+i}^{\ddag}$ satisfying
   \[
     |R_{\ell,k\ell+i}^{\ddag} | \leq     C (k-1)^{-2\alpha}\left(\varepsilon_\ell   + (k-1)^{-1}\right)
   \]
    for $\varepsilon_\ell \equiv   \sup_{j \geq \ell} |\gamma_{Z}(j)/[C_0 j^{-\alpha}] -1 | =o(1)$  and for a constant $C\geq 1$ not depending on $\ell$, $k \geq 2$, $0 \leq i \leq \ell-1$; the   bound on the remainder follows from (\ref{eqn:lcov}) (i.e., $|L_{\ell,m-1}| \leq C  $) along with $|\ell k + i + i_1 - i_2|\geq \ell(k-1)$ in the covariance bound $|\gamma_Z(\ell k + i + i_1 - i_2)|\leq C_0(1+\varepsilon_\ell) [\ell(k-1)]^{-\alpha}$ as well as $|\ell^{2\alpha}|\ell k + i + i_1 - i_2|^{-2\alpha} - k^{-2\alpha}|\leq 2 |k-1|^{-1-2\alpha}$
     for any $k \geq 2$, $0 \leq i, i_1,i_2 \leq \ell$ (noting $2\alpha \leq \am <1$ with $m \geq 2$).  Hence, we have
     \begin{eqnarray*}
     && s_{1,\ell}\\ &=& 4 C_0^{2}\left(  J_m^2 \frac{1}{(m-1)!}  L_{\ell,m-1} \right)^2
     \frac{\ell}{N}\sum_{k=2}^{\lceil n/\ell \rceil} (1-N^{-1}\ell k ) k^{-2\alpha} + O\left(  \frac{\ell}{N}+ \varepsilon_\ell  \frac{\ell}{N}\sum_{k=2}^{\lceil n/\ell \rceil}    k^{-2\alpha}  \right)  \\
      &=&  4 C_0^2  \left( J_m^2 \frac{1}{(m-1)!} \frac{2 C_0^{m-1}}{(1-\alpha (m-1))(2-\alpha_(m-1))}\right)^2  \frac{1}{(1-2\alpha)(2-2\alpha)}
      \left( \frac{\ell}{n}\right)^{2\alpha} (1+o(1))\\
      &\equiv& \phi_{\alpha,m}  \left( \frac{\ell}{n}\right)^{2\alpha} (1+o(1)),
     \end{eqnarray*}
     using $\varepsilon_\ell \rightarrow 0$ and that $2\alpha \leq \am <1$ and $0< (m-1)\alpha \leq \am <1$ (by $m \geq 2$) whereby
     \[
    L_{\ell,m-1}\rightarrow \frac{2 C_0^{m-1}}{(1-\alpha (m-1))(2-\alpha (m-1) )},
    \]
   and  $(\ell/N)^{1-2\alpha}\sum_{k=2}^{\lceil n/\ell \rceil} (1-N^{-1}\ell k )k^{-2\alpha}\rightarrow  [(1-2\alpha)(2-2\alpha)]^{-1}$  as $n\to \infty$ with $N/n\rightarrow 1$ and $\ell/n\rightarrow 0$.  For $m \geq 2$, this now establishes an   expansion for the first sum $s_{1,\ell}= \phi_{\alpha,m} (\ell/n)^{2\alpha}(1+o(1))$  in $t_{1,\ell} = s_{1,\ell}+s_{2,\ell}$.  Now (\ref{eqn:t3}) will follow for $t_{1,\ell}$ (when $m\geq 2$) by showing
    $s_{\ell,2} = o((\ell/n)^{2\alpha} )$.
    Considering the sum
    \[
   s_{\ell,2} \equiv  2  J_m^4 \sum_{A \in \mathcal{A}_c(m,m,m,m),\atop v_{12}\leq m-2} g(A)
     \frac{1}{N}\sum_{k=2\ell+1}^N  \left(1- \frac{k}{N}\right) R_{\ell,k}(A),
    \]
    note that, for a multigraph  $A \equiv (v_{12},v_{13},v_{14},v_{23},v_{24},v_{34})\in \mathcal{A}_c(m,m,m,m)$
    (with $v_{12}=v_{34}$, $v_{13}=v_{24}$, $v_{14}=v_{23}$, $v_{12}+v_{13}+v_{14}=m$), we may apply the bound from (\ref{eqn:RA2}) to find
    \begin{equation}
    \label{eqn:last}
 |R_{\ell, \ell k+i}(A)| \;\leq \; C  \frac{\ell^{2 \am} }{\ell^{ 2 \alpha v_{12} }}   [ C_1 \ell^{-\alpha} (k-1)^{-\alpha}]^{2m-2v_{12}}
 \;\leq \;  C C_1^{2m} (k-1)^{-4\alpha}
 \end{equation}
 holds for   generic constants $C>0, C_1 >1$, not depending on integers $\ell,k\geq 2$, $0\leq i \leq \ell-1$, and the values of $0 \leq v_{12} \leq m-2$; above we used $2m - 2v_{12} \geq 4$.  Application of the bound (\ref{eqn:last}) with $\sum_{A \in \mathcal{A}_c(m,m,m,m)} g(A)<\infty$ (e.g., $\mathcal{A}_c(m,m,m,m)$ is finite) then yields
 \[
   |s_{\ell,2}|  \;=\; O\left(\frac{\ell}{N}\sum_{k=1}^{\lceil n/\ell \rceil} k^{-4\alpha} \right) \;=\;
    \left\{ \begin{array}{lcl}
    O(\ell/n) && \mbox{if $\alpha >1/4$}\\
    O(\log(n/\ell)\cdot \ell/n) &&  \mbox{if $\alpha =1/4$}\\
       O( (\ell/n)^{4\alpha} ) &&  \mbox{if $\alpha <1/4$}
\end{array}     \right.\\
\; =\;
    o\left( \left(\frac{\ell}{N} \right)^{2\alpha}\right),
 \]
 which concludes the proof of (\ref{eqn:t3}) as well as the proof of Lemma~\ref{lem1}. $\Box$

\end{appendix}

\section*{Acknowledgements}
The authors are grateful to two anonymous referees and an Associate
Editor for constructive comments that  improved the
quality of this paper.
Research was  supported   by NSF DMS-1811998 and DMS-2015390.%
%


\begin{supplement}
\stitle{Proofs and other technical details}
\sdescription{A supplement \cite{ZLN22} contains  proofs and technical details along with further numerical results.}
\end{supplement}


\end{document}